\documentclass[aos,nameyear,dvips]{arximspdf}
\usepackage{dcolumn}
\usepackage{graphicx}

\doi{10.1214/09-AOS769}
\volume{38}
\issue{4}
\pubyear{2010}
\firstpage{2187}
\lastpage{2217}

\makeatletter

\newtheorem{theorem}{Theorem}
\newtheorem{corol}{Corollary}

\newproclaim{remark}{Remark}

\newtheorem{proposition}{Proposition}
\newtheorem{lemma}{Lemma}

\newcolumntype{k}[1]{D{,}{}{#1}}

\makeatother

\begin{document}
\begin{frontmatter}

\title{Nonparametric inference of quantile curves for nonstationary
time series\protect\thanksref{T1}}
\runtitle{Nonstationary quantile inference}

\thankstext{T1}{Supported in part by the NSF Grant DMS-04-78704.}

\begin{aug}
\author[A]{\fnms{Zhou} \snm{Zhou}\corref{}\ead[label=e1]{zhou@utstat.toronto.edu}}
\runauthor{Z. Zhou}
\affiliation{University of Toronto}
\address[A]{Department of Statistics\\
University of Toronto\\
100 St. George St.\\
Toronto, Ontario M5S 3G3\\
Canada\\
\printead{e1}} 
\end{aug}

\received{\smonth{8} \syear{2008}}
\revised{\smonth{10} \syear{2009}}

%
\begin{abstract}
The paper considers nonparametric specification tests of quantile
curves for a general class of nonstationary processes. Using
Bahadur representation and Gaussian approximation results for
nonstationary time series, simultaneous confidence bands and
integrated squared difference tests are proposed to test various
parametric forms of the quantile curves with asymptotically correct
type I error rates. A wild bootstrap procedure is implemented to
alleviate the problem of slow convergence of the asymptotic results.
In particular, our results can be used to test the trends of
extremes of climate variables, an important problem in understanding
climate change. Our methodology is applied to the analysis of the
maximum speed of tropical cyclone winds. It was found that an
inhomogeneous upward trend for cyclone wind speeds is pronounced at
high quantile values. However, there is no trend in the mean
lifetime-maximum wind speed. This example shows the effectiveness of
the quantile regression technique.
\end{abstract}

%

\end{frontmatter}

\section{Introduction}
\label{sec::intr}
After fitting a nonparametric model, one often asks whether it can
be simplified into certain parametric, semiparametric or more
parsimonious nonparametric forms. Recently, there has been an
enormous interest in developing nonparametric specification tests;
see, for example, Hall and Titterington (\citeyear{HT88}), Eubank and Speckman
(\citeyear{ES93}), H\"{a}rdle and Mammen (\citeyear{HM93}), Ingster
(\citeyear{I93}), Zheng (\citeyear{Z96}),
Hart (\citeyear{H97}), Stute (\citeyear{S97}), Xia (\citeyear{X98}),
Horowitz and Spokoiny (\citeyear{HS01}),
Fan, Zhang and Zhang (\citeyear{FZZ01}) and Fan and Jiang (\citeyear
{FJ07}) among others.
Many of the previous results concern nonparametric inference of the
(conditional) mean or density functions for independent data.

The primary goal of this paper is to perform nonparametric
specification tests of quantile curves for a class of nonstationary
processes that can be called locally stationary processes
[Draghicescu, Guillas and Wu (\citeyear{DGW08}) and Zhou and Wu
(\citeyear{ZW09a})].
Conceptually, the local stationarity is characterized by the
smoothly time-varying data generating mechanisms of the processes.
More precisely, let $\{X_{i,n}\}_{i=1}^n$ be the observed sequence.
We shall adopt the following formulation:
%
%
\begin{equation}\label{model}
X_{i,n} = G(i/n, \mathcal{F}_i),\qquad i = 1,2,\ldots,n,
\end{equation}
where $\mathcal{F}_i = (\ldots, \varepsilon_{i-1}, \varepsilon_i)$,
$\varepsilon_i$, $i \in\mathbb{Z}$, are independent and identically
distributed (i.i.d.) random variables, and $G\dvtx[0,1] \times\mathbb{R}
^\infty
\mapsto\mathbb{R}$ is a measurable function such that $\zeta_i(t) := G(t,
\mathcal{F}_i)$ is a properly defined random variable for all $t\in
[0,1]$. Here, local stationarity means that the random function
$\zeta_i(t)$ is smooth in $t \in[0, 1]$ in an appropriate sense.
Process (\ref{model}) covers a wide range of nonstationary linear
and nonstationary nonlinear processes and it naturally extends many
existing stationary time series models into the nonstationary
setting. See Zhou and Wu (\citeyear{ZW09a}) for more discussion and examples.
In the sequel, for notational convenience, we shall write $X_{i,n}$
as $X_i$.

The paper is motivated by the important problem of understanding the
trends of extremes and variability of climate variables. As stated
in Katz and Brown (\citeyear{KB92}), ``\textit{understanding climate change
demands attention to changes in climate variability and extremes}.''
By interpreting climate extremes as upper and lower quantiles and
climate variability as interpercentile ranges, the nonparametric
quantile estimation [Koenker (\citeyear{K05})] provides a simple and effective
means to address the latter problem. On the other hand,
climatologists often want to know whether a simple linear or
quadratic function is appropriate to describe the trends. The
parametric description of the trends is more interpretable and
efficient than its nonparametric counterparts if the parametric
model is correctly specified. To address the latter issue, it is
necessary to develop nonparametric specification tests for the
quantile curves.




For independent data, there have been a few results on nonparametric
specification tests of quantile curves. He and Zhu (\citeyear{HZ03})
and Kim
(\citeyear{K07}) proposed tests based on the cusum process of the residuals
and the Rao-score statistics, respectively. The tests use fitted
values and residuals under the null hypothesis (namely the
parametric model) and therefore enjoy the advantage of avoiding the
need for nonparametric function estimation. On the other hand,
however, the above tests are sensitive only to a restrictively small
class of alternatives; namely alternatives in the neighborhood of
the parametric null. An alternative specification test which
achieves minimax rate over a large class of smooth functions is
proposed in Horowitz and Spokoiny (\citeyear{HS02}). It seems that the
test is
tailored to the specific problem of testing linearity. For other
contributions, see Zheng (\citeyear{Z98}), Rosenkrantz (\citeyear
{R00}) and Wang (\citeyear{W07},
\citeyear{W08}).

To our knowledge, testing parametric functional forms of quantile
curves under the time series setting has not yet been considered in
the literature. In this paper, we shall conduct the specification
tests by directly comparing the nonparametric quantile estimates and
the parametric null. In principle, if some global measure such as
the $\mathcal{L}^{\infty}$ norm or the $\mathcal{L}^2$ norm is large for
the difference between the nonparametric fits and the parametric
null, then there is evidence against the parametric null hypothesis.
This approach is intuitively plausible and it is easy to understand.
For nonparametric inference of the mean functions, this very idea
has been widely applied; see, for example, Eubank and Speckman
(\citeyear{ES93})
and H\"{a}rdle and Mammen (\citeyear{HM93}) among others.

Nevertheless, obtaining the asymptotic distributions of
appropriately normalized global measures of the deviations of the
quantile curves has been a very difficult problem when dependence is
present. Generally speaking, the latter problem can be solved if we
have (i) unform Bahadur representations of the estimated quantile
curves and (ii) a sharp Gaussian approximation result for the
weighted empirical processes of the nonstationary time series
$(X_i)_{i=1}^n$. Recently, Zhou and Wu (\citeyear{ZW09a}) obtained a sharp
uniform Bahadur representation for the local linear quantile
estimates of locally stationary time series. On the other hand, Wu
and Zhou (\citeyear{WZ09}) established Gaussian approximation results for
partial sums of nonstationary time series with nearly optimal
rates. An extension of the latter result addresses the above issue
(ii); see Theorem \ref{thm::sip} below.

With the recent progress in Bahadur representations and Gaussian
approximations for locally stationary time series, we are able to
construct in this paper simultaneous confidence bands (SCB) (based
on the $\mathcal{L}^{\infty}$ norm) and integrated squared difference
tests (ISDT) (based on the $\mathcal{L}^2$ norm) for the quantile
curves as tools for the nonparametric inference. The SCB is shown to
asymptotically achieve the correct coverage probability. Moreover,
we prove that the ISDT asymptotically attain the correct type I
error rate and are asymptotically optimal in terms of rates of
convergence for nonparametric hypothesis testing in the sense of
Ingster (\citeyear{I93}). Our results shed new light on nonparametric
specification tests of $M$-type estimates for nonstationary time
series.

The rest of the paper is structured as follows. Section
\ref{sec::pre} introduces the local linear quantile estimates and
the dependence measures. Section \ref{sec::mr} presents the
asymptotic results on the nonparametric specification tests. In
particular, asymptotic results on the SCB and ISDT are presented in
Sections \ref{sec::scb} and \ref{sec::l2}, respectively. A wild
bootstrap procedure is introduced in Section
\ref{sec::implementation} to address the issue of slow convergence
of the SCB and ISDT tests. Section \ref{sec::implementation} also
contains discussions on bandwidth selection and nuisance parameter
estimation. Section \ref{sec::simu} presents two simulation studies
to compare accuracy and sensitivity of various nonparametric tests.
The global tropical cyclone data is studied in Section
\ref{sec::data}. Proofs are given in Section \ref{sec::pr}.

\section{Preliminaries}
\label{sec::pre} We now introduce some notation. For a vector
$\mathbf{v} = (v_1, v_2, \ldots, v_p) \in\mathbb{R}^p$, let
$|\mathbf{v}| =
(\sum_{i=1}^p v_i^2)^{1/2}$. For a $p\times p$ matrix $A$, define
$|A| = \sup\{ |A\mathbf{v}|\dvtx|\mathbf{v}|=1\}$. For a random vector
$\mathbf{V}$, write $\mathbf{V} \in\mathcal{L}^q$ ($q > 0$) if $\|
\mathbf{V}\|_q := [\mathbb{E}(|\mathbf{V}|^q) ]^{1/q} < \infty$ and
$\|
\mathbf{V}\| =
\|\mathbf{V}\|_2$. Denote by $\Rightarrow$ the weak convergence. For an
interval $\mathcal{I}\subset\mathbb{R}$, denote by $\mathcal
{C}^i\mathcal{I}$,
$i \in
\mathbb{N}$, the collection of functions that have $i$th order continuous
derivatives on $\mathcal{I}$, and, for $\mathcal{D} \subset\mathbb
{R}^d$, let
$\mathcal{C}\mathcal{D}$ be the collection of real-valued functions that
are continuous on $\mathcal{D}$. A function $f\dvtx\mathbb
{R}^d\rightarrow
\mathbb{R}$ is
said Lipschitz continuous on $\mathcal{D} \subset\mathbb{R}^d$ if
there exists
a finite constant $C$, such that $|f(x_1)-f(x_2)| \le C |x_1-x_2|$
for all $x_1, x_2 \in\mathcal{D}$. For $x \in\mathbb{R}$, define
$x^+ =
\max(x,0)$. The symbol $C$ denotes a finite generic constant which
may vary from line to line.

\subsection{Local linear quantile estimator}\label{sec::quantile_estimate}
Recall $\zeta_i(t) = G(t,
\mathcal{F}_i)$. Let $F(t,x) = \mathbb{P}(\zeta_i(t) \le
x)$, $x \in\mathbb{R}$, be the cumulative distribution function (cdf) of
$\zeta_i(t)$, $t \in[0, 1]$; let $Q_{\alpha}(t)$ be the $\alpha$th
quantile function of $\zeta_i(t)$, $\alpha\in(0, 1)$, namely
$Q_\alpha(t) = \inf_x \{F(t,x) \ge\alpha\}$. Suppose $Q_{\alpha}(t)$
is smooth on $[0,1]$. As $Q_{\alpha}(t_1) \approx
Q_{\alpha}(t) + (t_1-t)Q_{\alpha}'(t)$ for $t_1$ close to $t$,
$Q_\alpha
(t)$ and $Q'_\alpha(t)$ can be estimated by the local linear approach
[Koenker (\citeyear{K05})]
\[
(\hat{Q}_{\alpha,b_n}(t),\hat{Q}'_{\alpha,b_n}(t))=\mathop{\arg
\min}
_{(\beta_0,\beta_1)}
\sum_{i=1}^n\rho_{\alpha}\bigl(X_i-\beta_0-\beta_1(t_i-t)\bigr)
K_{b_n}(t_i-t),
\]
where $t_i=i/n$, $\rho_{\alpha}(x) = \alpha x^++(1-\alpha)(-x)^+$ is
the check
function. Here $K$ is a kernel function, $K_{b_n}(\cdot) = K(\cdot/ b_n)$
and $b_n=b_n(\alpha) > 0$ is the bandwidth depending on $\alpha$. We
shall omit the
subscript $b_n$ in $\hat{Q}$ and $\hat{Q}'$ hereafter if no
confusion will be caused.

\subsection{The dependence measures}\label{sec::model}

It is shown [Condition (B1) of Zhou and Wu (\citeyear{ZW09a})] that the
asymptotic behavior of the local linear quantile estimator is
determined by the dependence structure of
\[
F_k(t,x,\mathcal{F}_i):=\frac{\partial^k\mathbb{P}\{G(t,\mathcal
{F}_{i+1})\le
x|\mathcal{F}_i\}}{\partial x^k},\qquad k=0,1,2,3.
\]
To quantify the dependence structure of the above processes,
let us consider a generic nonlinear system $\{H(t,x,\mathcal{F}_i)\}
_{i\in\mathbb{Z}}$, where
$(t,x)\in[0,1] \times\mathbb{R}$ and $H\dvtx[0,1] \times\mathbb
{R}\times\mathbb{R}
^\infty
\mapsto\mathbb{R}$ is a measurable function such that $H(t,x,\mathcal
{F}_i)$ is
well defined for all $(t,x) \in[0,1]\times\mathbb{R}$.
Let $\varepsilon'_i$, $i \in\mathbb{Z}$, be an i.i.d. copy of
$\varepsilon_j$,
$j \in\mathbb{Z}$.
For $k \ge0$, let $\mathcal{F}^*_k = (\mathcal{F}_{-1}, \varepsilon'_0,
\varepsilon_1, \ldots, \varepsilon_k)$ and define the physical
dependence measure
%
%
\begin{equation}\label{eq::ph-de}
\delta_H(k,p)={\sup_{(t,x)\in[0,1]\times\mathbb{R}}}
\|H(t,x,\mathcal{F}_k)-H(t,x,\mathcal{F}^*_{k})\|_p.
\end{equation}
Here we recall $\|\cdot\|_p=[\mathbb{E}(|\cdot|^p) ]^{1/p}$. Note that
$\delta_H(k,p)$ measures the overall dependence of $H(t,x, \mathcal
{F}_k)$ on the input $\varepsilon_0$. The physical dependence
measures are by their definition closely related to the data
generating mechanism and hence are easy to work with; see Section 4
of Zhou and Wu (\citeyear{ZW09a}) for the related calculations for locally
linear and nonlinear time series models.

We shall call the system $\{H(\cdot,\cdot,\mathcal{F}_i)\}_{i\in
\mathbb{Z}}$
uniformly geometric moment contracting of order $p$ [UGMC($p$)] if
$\delta_H(k,p)$ decays exponentially with respect to $k$; namely
%
%
\begin{equation}\label{eq::gmc}
\delta_H(k,p)=O(\chi^k),\qquad 0<\chi<1.
\end{equation}
Following Section 4 of Zhou and Wu (\citeyear{ZW09a}), condition
(\ref{eq::gmc}) is readily verifiable for a large class of
nonstationary nonlinear processes and nonstationary linear models.
All our results will be presented in terms of the physical
dependence measures and the UGMC conditions.

Note that UGMC($2$) is stronger than the stability condition of Zhou
and Wu (\citeyear{ZW09a}). Analogous results of this paper can be
proved with
$\delta_H(k,p)$ decaying algebraically at a sufficiently fast rate.
However, the technical details are much lengthier and we chose to
use the UGMC condition for clarity of presentation.

\section{Main results}\label{sec::mr}

\subsection{Assumptions}\label{sec::sip}
We shall make the following assumptions on the process $(X_i)$ and the
kernel $K$.
\begin{enumerate}[(A2)]

\item[(A1)] $f(t,x)$ is Lipschitz continuous on $[0,1]\times\mathbb
{R}$, where
$f(t, \cdot)$ is the density function of $\zeta_i(t)$. Assume $\inf
_{t\in[0,1]}f(t,Q_{\alpha}(t))>0$.

\item[(A2)] Let $J(t, x, \mathcal{F}_i) = I\{G(t,\mathcal{F}_i)\le
x\}$ and
%
%
\begin{equation}\label{eq::lrv}
\sigma^2(t)=\sum_{i=-\infty}^{\infty}\operatorname{cov}[J(t,Q_{\alpha
}(t),\mathcal{F}_0),J(t,Q_{\alpha}(t),\mathcal{F}_i)].
\end{equation}
Assume $\sigma(t)$ is Lipschitz continuous on $[0,1]$ and $\inf_{t\in
[0,1]}\sigma(t)>0$.

\item[(A3)] (Stochastic Lipschitz continuity condition.) There exists
$q\ge1$, such that $\| \zeta_i(t_1) - \zeta_i(t_2) \|_q \le C
|t_1-t_2|$ holds for all $t_1, t_2 \in[0,1]$.

\item[(A4)]
Assume that, for $k=1,2,3$, $F_k(t,x,\mathcal{F}_i) := \partial^k F(t,
x, \mathcal{F}_i)/\partial x^k$ exists and
$\sum_{k=0}^{\infty}\delta_{F_i}(k,4)<\infty$, for $i=1,2,3$.

\item[(A5)] $F(t, x, \mathcal{F}_i)$ is UGMC(4).

\item[(A6)] There exists $C_0 < \infty$ such that
$\sup_{(t,x) \in[0,1] \times\mathbb{R}} F_1(t,x,\mathcal{F}_i)<
C_0$ almost
surely.

\item[(K1)] $K(\cdot) \in\mathcal{K}$, where $\mathcal{K}$ is the
collection of density functions $K$ such that $K$ is
symmetric with support $[-1, 1]$ and $K \in\mathcal{C}^1[-1,1]$.
For $K(\cdot) \in\mathcal{K}$ and $j\ge1$, define
\begin{eqnarray*}
\phi_{K} &=& \int_{-1}^{1} K^2(x) \,dx, \qquad \mathfrak{C}_K=\frac{{\int
_{-1}^1}|K'(x)|^2\,dx}{\phi_K}, \\
\mu_{j,K} &=& \int_{-1}^{1} x^j K(x) \,dx.
\end{eqnarray*}
We shall omit the subscript $K$ in the sequel if no confusion will be caused.
\end{enumerate}

A few remarks on the above regularity conditions are in order. The
function $\sigma^2(t)$ in condition (A2) is called the long-run
variance function of $\{J(t,Q_{\alpha}(t)$, $\mathcal{F}_i)\}_{i=-\infty
}^{\infty}$ to account for the dependence of the
series. Note that $\sigma^2(t)$ is well defined even for
heavy-tailed processes $(X_i)$ since $0\le I\{G(t,\mathcal{F}_i)\le
x\}\le1$. Condition (A3) means local stationarity and it asserts
smoothness of $G(t, \mathcal{F}_i)$ with respect to time~$t$.
Conditions (A4) and (A5) assert that processes $\{F_k(t,x,\mathcal
{F}_i)\}$, $k=0,1,2,3$, are short range dependent (SRD). Conditions
(A3)--(A5) can be verified for a large class of nonstationary linear
and nonlinear processes by the arguments in Section 4 of Zhou and Wu
(\citeyear{ZW09a}). Condition (A6) is mild and it means that the conditional
density function of $J(t, x, \mathcal{F}_i)$ is bounded. A popular
choice of the kernel function is the Epanechnikov kernel $K(x) =
3\max(0, 1-x^2) /4$.

\subsection{Simultaneous confidence bands}\label{sec::scb}
The simultaneous confidence band\break (SCB) is a classic tool for
nonparametric inference. To construct a $100(1-\beta)\%$ SCB for
$Q_{\alpha}(\cdot)$, one finds two functions $l$ and $u$ depending
on $(X_i)_{i=1}^n$, such that
\[
\lim_{n\rightarrow\infty}\mathbb{P}\bigl(l(t)\le Q_{\alpha}(t)\le
u(t)\mbox{
for all
} t\in(0,1)\bigr)=1-\beta.
\]
A candidate function for $Q_{\alpha}(\cdot)$ is rejected at level
$\beta$ if it is not fully contained in the SCB. The SCB provides
appreciable direct visual information on the overall variability of the
fitted curves. See, for example, Bickel and Rosenblatt (\citeyear
{BR73}) for the inference of density functions; Eubank and Speckman
(\citeyear {ES93}), Sun and Loader (\citeyear{SL94}), Neumann and
Kreiss (\citeyear{NK98}), Wu and Zhao (\citeyear{WZ07}) and Zhao and Wu
(\citeyear{ZW08}) for the inference of (conditional) mean functions and
Fan and Zhang (\citeyear{FZ00}) for the inference of coefficient
functions of varying coefficient models. In this section, we shall
establish the asymptotic theory for the maximal absolute deviation of
$\hat{Q}_{\alpha}(t)$ from $Q_{\alpha}(t)$ on $(0,1)$. The theoretical
results facilitate construction of a SCB of $Q_{\alpha}(t)$ which
asymptotically achieves the nominal coverage probability.
\begin{theorem}\label{thm::scb}
Assume $Q_{\alpha}(\cdot)\in\mathcal{C}^3[0,1]$ and conditions
\textup{(A1)--(A6)} and \textup{(K1)} hold. Further assume $\sqrt{n}b_n/\log^5
n\rightarrow\infty$ and $nb_n^7\log n\rightarrow0$, then we have
%
%
\begin{eqnarray}\label{eq::scb}
&&\lim_{n\rightarrow\infty}\mathbb{P}\biggl[\sup_{t\in\mathcal{T}_n} \biggl\{
\frac
{\sqrt{nb_n}f(t,Q_{\alpha}(t))}{\sqrt{\phi}\sigma(t)}\nonumber\\
&&\hspace*{56.9pt}{}\times|\hat
{Q}_{\alpha
}(t)-Q_{\alpha}(t)-\mu_2b_n^2Q''_{\alpha}(t)/2| \biggr\}\\
&&\hspace*{56.9pt}\hspace*{62.6pt}{}-B(m^*)\le\frac{x}{\sqrt{2\log m^*}}
\biggr]= e^{-2e^{-x}},\nonumber
\end{eqnarray}
where $\mathcal{T}_n=[b_n,1-b_n]$, $m^*=1/b_n$ and
\[
B(m^*)=(2\log m^*)^{1/2}+(2\log m^*)^{-1/2}[\log\mathfrak{C}-2\log
\pi
-2\log2]/2.
\]
\end{theorem}

For a fixed level $\beta$, Theorem \ref{thm::scb} implies that one
can construct a $100(1-\beta)\%$ simultaneous confidence band for
$Q_{\alpha}(t)$
%
%
\begin{equation}\label{eq::asySCB}
\hat{Q}_{\alpha}(t)-\mu_2b_n^2Q''_{\alpha}(t)/2 \pm\frac{\sqrt
{\phi
}\sigma(t)}{\sqrt{nb_n}f(t,Q_{\alpha}(t))}\biggl(B_K(m^*)+\frac{u_{\beta
}}{\sqrt{2\log m^*}}\biggr),
\end{equation}
where $u_{\beta}=-\log\log[(1-\beta)^{-1/2}]$. By Theorem
\ref{thm::scb}, SCB (\ref{eq::asySCB}) asymptotically achieves the
right coverage probability $1-\beta$. Note that
$V(t):=\frac{\sqrt{\phi}\sigma(t)}{\sqrt{nb_n}f(t,Q_{\alpha}(t))}$
is the asymptotic standard deviation of $\hat{Q}_{\alpha}(t)$
[Theorem 1 of Zhou and Wu (\citeyear{ZW09a})].

The following theorem concerns the local power of the SCB (\ref{eq::asySCB}).
\begin{theorem}\label{thm::scb_power}
Suppose $Q_{\alpha}(t)=Q^o_{\alpha}(t)+\gamma_n\eta(t)+o(\gamma_n)$,
where $Q^o_{\alpha}(t), \eta(t)\in\mathcal{C}[0,1]$,
$\gamma_n=1/\sqrt{-2nb_n\log b_n}$ and $o(\gamma_n)$ is uniform in
$t$ on $[0,1]$. Then under conditions of Theorem \ref{thm::scb}, we
have
%
%
\begin{eqnarray}\label{eq::scb_power}
&&\lim_{n\rightarrow\infty}\mathbb{P}\biggl[\sup_{t\in\mathcal{T}_n} \biggl\{
\frac
{\sqrt{nb_n}f(t,Q_{\alpha}(t))}{\sqrt{\phi}\sigma(t)}\nonumber\\
&&\hspace*{56.9pt}{}\times|\hat
{Q}_{\alpha
}(t)-Q^o_{\alpha}(t)-\mu_2b_n^2Q''_{\alpha}(t)/2| \biggr\}\\
&&\hspace*{56.9pt}\hspace*{62.6pt}{}-B(m^*)\le\frac{x}{\sqrt{2\log m^*}} \biggr]=
e^{-s(\eta)e^{-x}},\nonumber
\end{eqnarray}
where
\[
s(\eta)=\int_0^1\exp\biggl\{\frac{\eta(t)f(t,Q_{\alpha}(t))}{\sqrt
{\phi
}\sigma(t)} \biggr\} \,d t+\int_0^1\exp\biggl\{\frac{-\eta(t)f(t,Q_{\alpha
}(t))}{\sqrt{\phi}\sigma(t)} \biggr\} \,dt.
\]
\end{theorem}

Theorem \ref{thm::scb_power} follows from similar arguments to those
in the proof of Theorem \ref{thm::scb} and Theorem A1 of Bickel and
Rosenblatt (\citeyear{BR73}). Details are omitted.

Theorem \ref{thm::scb_power} implies that SCB (\ref{eq::asySCB}) can
detect alternatives with the rate $\gamma_n$. For the
$100(1-\beta)\%$ SCB (\ref{eq::asySCB}), the asymptotic power of the
test $H_0\dvtx Q_{\alpha}(t)=Q^o_{\alpha}(t)$ versus $H_a\dvtx
Q_{\alpha}(t)\neq Q^o_{\alpha}(t)$ is
%
%
\begin{equation}\label{eq::power}
1-(1-\beta)^{s(\eta)/2}
\end{equation}
under the local alternatives specified in Theorem
\ref{thm::scb_power}. Since $s(\eta)\ge2$ and $s(\eta)=2$ if and
only if $\eta(t)\equiv0$, our test based on the SCB is always
asymptotically unbiased under such alternatives. In the case of
density function inference for i.i.d. data, the same result was
obtained by Bickel and Rosenblatt (\citeyear{BR73}).

\subsubsection{Optimality of the SCB}\label{sec::op_scb}
If the local linear smoothing technique is adopted and the bandwidth
series $(b_n)$ is fixed, then SCB (\ref{eq::asySCB}) is optimal in
the sense that asymptotically it covers the minimum area. To see
this, a Lagrange multiplier argument can be implemented. A similar
argument can be found in Zhou and Wu (\citeyear{ZW09b}) for nonparametric
inference of time-varying coefficients in functional linear models.
Note by equations (\ref{eq::thm115}) and (\ref{eq::thm116}) in
Section \ref{sec::pr}, we have under conditions of Theorem
\ref{thm::scb},
%
%
\begin{equation}\label{eq:scb_approx}\quad
\sup_{t\in\mathcal{T}_n}|[\hat{Q}_{\alpha}(t)-Q_{\alpha
}(t)-\mu_2b_n^2Q''_{\alpha}(t)/2]/V(t)- \Theta_n(t)|=o_\mathbb
{P}(\log
^{-1/2} n),
\end{equation}
where $\Theta_n(t)=\sum_{i=1}^nV_iK_{b_n}(t_i-t)/\sqrt{\phi nb_n}$
with $(V_i)_{i=1}^n$ i.i.d. standard normal. In other words, the
simultaneous fluctuations of
$[\hat{Q}_{\alpha}(t)-Q_{\alpha}(t)-\mu_2b_n^2Q''_{\alpha}(t)/2]/V(t)$
are asymptotically equivalent to those of $\Theta_n(t)$. Let $s_i =
(2i-1)b_n$, $i=1, \ldots, g_n$, where $g_n = \lfloor1 / (2b_n)
\rfloor$. From (\ref{eq:scb_approx}),
$\hat{Q}_{\alpha}(s_i)-Q_{\alpha}(s_i)-\mu_2b_n^2Q''_{\alpha}(s_i)/2$
are asymptotically independent $N(0,V^2(s_i))$. Suppose a band
\[
l(s_i)\le\hat{Q}_{\alpha}(s_i)-Q_{\alpha}(s_i)-\mu
_2b_n^2Q''_{\alpha
}(s_i)/2\le u(s_i),\qquad i=1,2,\ldots,g_n,
\]
achieves a preassigned coverage probability $1-\beta$. From the
above discussion, we see that the coverage probability restriction
can be asymptotically written as
\[
c(n,b_n):=\prod_{i=1}^{g_n} \bigl[\Phi\bigl(u(s_i)/V(s_i)\bigr)
-\Phi\bigl(l(s_i)/V(s_i)\bigr) \bigr]=1-\beta,
\]
where $\Phi(\cdot)$ is the normal cumulative distribution function
(cdf). In order to achieve the minimum average length, one minimizes
the following Lagrange multiplier:
%
%
\begin{equation}
\sum_{i=1}^{g_n}[u(s_i)-l(s_i)]-\lambda\{\log[c(n,b_n)]-\log
(1-\beta)\}.
\end{equation}
Simple calculations show that the minimum is achieved at
$u(s_i)=-l(s_i)=g(n,b_n,\beta)V(s_i)$, where $g(n,b_n,\beta)$ is a
deterministic function. The important message here is that the
asymptotically optimal SCB at each time point $t$ should have length
proportional to the asymptotic standard deviation of
$\hat{Q}_{\alpha}(t)$, which is the case in our construction.

\subsection{The integrated squared difference test (ISDT)}\label{sec::l2}
Another popular basis for nonparametric inference is $\mathcal
{L}^2$-distance based tests. In general, one calculates a $\mathcal
{L}^2$ norm related distance between the fitted nonparametric curve
and the parametric null, and a large distance indicates violation of
the null hypothesis. Most of the existing results on the $\mathcal
{L}^2$ type tests are for independent data. See, for instance, Bickel
and Rosenblatt (\citeyear{BR73}), H\"{a}rdle and Mammen (\citeyear
{HM93}), Zheng (\citeyear{Z96}),
Fan, Zhang and Zhang (\citeyear{FZZ01}), Zhang and Dette (\citeyear
{ZD04}) and Fan and
Jiang (\citeyear{FJ05}) among others.


A simple way to construct a $\mathcal{L}^2$ type test is to use the
statistic
$T_n=\int_0^1[\hat{Q}_{\alpha}(t)-Q_{\alpha}(t)]^2\pi(t) \,dt$.
However, the bias of the local linear estimate $\hat{Q}_{\alpha}(t)$
is of order $O(b_n^2+\frac{1}{nb_n})$ and is not negligible for the
asymptotic analysis. The extra bias term complicates the asymptotic
distribution and reduces the precision of the test $T_n$. See also
H\"{a}rdle and Mammen (\citeyear{HM93}) for a related discussion in
the case of
conditional mean inference. To overcome the disadvantage, we shall
use the following jackknife bias reduction technique [Wu and Zhao
(\citeyear{WZ07})]:
%
%
\begin{equation}\label{eq::jackknife}
\tilde{Q}_{\alpha,b_n}(t)=2\hat{Q}_{\alpha,b_n}(t)-\hat{Q}_{\alpha
,\sqrt
{2}b_n}(t).
\end{equation}
It can be shown that bias of $\tilde{Q}_{\alpha}(t)$ is of order
$o(b_n^2+\frac{1}{nb_n})$ uniformly on $[0,1]$ if
$Q_{\alpha}(t)\in\mathcal{C}^2[0,1]$. Note that using
(\ref{eq::jackknife}) is equivalent to using the second-order kernel
\[
K^*(x):=2K(x)-K\bigl(x/\sqrt{2}\bigr)/\sqrt{2}.
\]
To test $H_0\dvtx Q_{\alpha}(t)=Q_{\alpha}^o(t)$ versus $H_a\dvtx
Q_{\alpha}(t)\neq Q_{\alpha}^o(t)$, we propose the following test
statistic
%
%
\begin{equation}\label{eq::L2test}
T_n^*=\int_{\mathcal{T}_n^*}[\tilde{Q}_{\alpha}(t)-Q^o_{\alpha
}(t)]^2\pi
(t) \,dt,
\end{equation}
where $\mathcal{T}_n^*=[\sqrt{2}b_n,1-\sqrt{2}b_n]$ and the weight
$\pi(t)$ are assumed to be nonnegative and Lipschitz continuous in
$t\in[0,1]$. The following theorem establishes the asymptotic
normality of $T_n^*$.
\begin{theorem}\label{thm::L2test}
Assume that $Q_{\alpha}(t)\in\mathcal{C}^2[0,1]$; that conditions
\textup{(A1)--(A6)} and \textup{(K1)} hold and that
$Q_{\alpha}(t)=Q^o_{\alpha}(t)+\varrho_n\eta(t)+o(\varrho_n)$, where
$Q^o_{\alpha}(t), \eta(t)\in\mathcal{C}[0,1]$,
$\varrho_n=n^{-1/2}b_n^{-1/4}$ and $o(\varrho_n)$ is uniform in $t$
on $[0,1]$. Further assume $nb_n^4/\log^{10} n\rightarrow\infty$ and
$nb_n^{9/2}=O(1)$. We have
%
%
\begin{eqnarray}\label{L2normality}
&&n\sqrt{b_n}T_n^*-\frac{1}{\sqrt{b_n}}K^*\star K^*(0)\int_{0}^1\pi
^*(t) \,dt-\int_0^1\eta^2(t)\pi(t) \,dt\nonumber\\[-8pt]\\[-8pt]
&&\qquad\Rightarrow N \biggl(0,2\int_\mathbb{R}[K^*\star K^*(t)]^2 \,dt\int
_{0}^1\pi
^*(t)^2 \,dt \biggr),\nonumber
\end{eqnarray}
where $\star$ denotes the convolution operator and $\pi^*(t)=\pi
(t)\sigma^2(t)/f^2(t,Q_{\alpha}(t))$.
\end{theorem}

When $\eta(t)\equiv0$, Theorem \ref{thm::L2test} unveils the
asymptotic null distribution of $T_n^*$. The ISDT can detect
alternatives with the rate $\varrho_n$. Under the local alternatives
specified in Theorem \ref{thm::L2test}, simple calculations based on
(\ref{L2normality}) show that the asymptotic power of the ISDT with
level $\beta$ equals
%
%
\begin{equation}\label{eq::powerl2}
\Phi\biggl(\frac{\int_0^1\eta^2(t)\pi(t) \,dt}{[2\int_\mathbb
{R}[K^*\star
K^*(t)]^2 \,dt\int_{0}^1\pi^*(t)^2 \,dt]^{1/2}}-z_{1-\beta} \biggr),
\end{equation}
where $\Phi(\cdot)$ and $z_{1-\beta}$ denote the cumulative
distribution function and the $1-\beta$ quantile of the standard
normal distribution. Therefore, a simple use of the Cauchy--Schwarz
inequality shows that choosing weights $\pi(t)$ proportional to
$\eta^2(t)f^4(t,Q_{\alpha}(t))/\sigma^4(t)$ maximizes the above
asymptotic power. Of course, in real applications it is difficult to
specify $\eta(t)$. Therefore, one can simply choose $\pi(t)\equiv1$.
On the other hand, we suggest choosing
$\pi(t)=f^2(t,\hat{Q}_{\alpha}(t))/\hat{\sigma}^2(t)$, which leads
to an easier implementation of the bootstrap. See the discussions in
Section \ref{sec::boots} for more details.
\begin{corol}\label{cor::optimality}
Under conditions\vspace*{2pt} of Theorem \ref{thm::L2test}, $T_n^*$ can detect
alternatives with departure rate $n^{-4/9}$ if the bandwidth
$b_n=O(n^{-2/9})$.
\end{corol}
\begin{remark}\label{rem::optimality}
If $X_i$ is distributed as $N(\mu(t_i),\sigma^2)$, and the
$X_i$'s are independent, then inference of quantile curves is
equivalent to inference of the mean function~$\mu(\cdot)$. In this
case, Ingster (\citeyear{I93}) and Lepski and Spokoiny (\citeyear
{LS99}) proved that the
optimal rate for testing $H_0\dvtx Q_{\alpha}(t)=Q_{\alpha}^o(t)$ is
$n^{-4/9}$. Hence, the integrated squared difference test $T_n^*$ is
optimal in the sense that it achieves the optimal rate of
convergence. Furthermore, Corollary \ref{cor::optimality} implies
that for nonparametric quantile function testing, weak dependence
and local stationarity do not deteriorate the optimal rate. However,
it should be noted that when there is long memory, the rate will be
deteriorated.
\end{remark}
\begin{remark}\label{rem::1}
Since $\gamma_n \gg\varrho_n$, the ISDT test $T_n^*$ dominates
the test based on SCB (\ref{eq::asySCB}) if bandwidths of the same
order are used for the tests; namely $T_n^*$ is asymptotically more
powerful. Therefore a general rule of thumb is to use the SCB when
one wants to explore the overall pattern of the quantile curves and
to implement the ISDT test when one is interested in verifying a
specific parametric null. On the other hand, as stated in H\"{a}rdle
and Mammen (\citeyear{HM93}), ``\textit{certainly from a more data
analytic point
of view distances would be more satisfactory which reflect
similarities in the shape of the regression functions}.'' For a
moderate sample size, intuitively, the ISDT test would be more
powerful compared to the SCB test if the true quantile curve differs
from the null in a systematic and even way; while the SCB test is
better when the latter difference is abrupt or ``bumpy,'' in which
case the $\mathcal{L}^2$ norm does not reflect the characteristics of
the difference. Hence, if there is some prior knowledge on the shape
of the discrepancy, one could select a test accordingly. In Section
\ref{sec::simu}, we shall conduct a simulation study to compare the
powers of the two tests under various alternatives.
\end{remark}

\section{Implementation}\label{sec::implementation}

\subsection{A wild bootstrap procedure}\label{sec::boots}
The asymptotic results in Section \ref{sec::mr} are based on the
uniform Gaussian approximations of
$\hat{Q}_{\alpha}(t)-Q_{\alpha}(t)$ on $(0,1)$. It is known that the
convergence rates of the $\mathcal{L}^{\infty}$ and $\mathcal{L}^2$ norms
of the corresponding Gaussian processes are very slow (see also
proofs in Section \ref{sec::pr}). For example, when
$b_n(\alpha)=O(n^{-1/5})$, the convergence rates of the $\mathcal
{L}^{\infty}$ and $\mathcal{L}^2$ norms are $1/\log^{1/2} n $ and
$n^{-1/10}$, respectively. Therefore, for moderate sample sizes, tests
based on the asymptotic theory are not reliable.

For nonparametric inference, the bootstrap is a classic tool for
achieving faster convergence. It is impossible to have a complete
list of literature here and we shall only mention several
representatives. Among others, Mammen (\citeyear{M93}), H\"{a}rdle and Mammen
(\citeyear{HM93}), Chapter 8 of Shao and Tu (\citeyear{ST95}),
Stute, Gonzalez Manteiga and Presedo Quindimil (\citeyear{SGmPq98}),
Neumann and Kreiss (\citeyear{NK98}) and Fan and Jiang (\citeyear
{FJ07}) considered wild
bootstrap inference for conditional mean regression for independent
data; Jhun (\citeyear{J88}), Faraway and Jhun (\citeyear{FJ90}) and
Hall (\citeyear{H93})
considered bootstrap inference of the density function for i.i.d. data.
For nonparametric inference of dependent data, among others, Politis
and Romano (\citeyear{PR94}) proposed a stationary bootstrap for simultaneous
inference of the spectral density functions of weakly dependent
stationary time series and Wu and Zhao (\citeyear{WZ07}) used a wild bootstrap
technique to test the mean function under stationary errors. For
other contributions, see Barrio and Matr\'{a}n (\citeyear{BM00}).

On the other hand, there is also large literature on bootstrap
methods for parametric quantile regression. See, for example, Chapter
3.9 of Koenker (\citeyear{K05}) and references therein for independent data
and Fitzenberger (\citeyear{F98}) for the moving block bootstrap for strong
mixing samples.


Despite the huge literature on bootstrap strategies for
nonparametric inference and parametric quantile regression, there
have been few results on bootstrap methods for nonparametric
quantile inference for nonstationary time series. The major
difficulty lies in the fact that joint distributions of subseries
within different time spans can be drastically different for a
nonstationary process; therefore it is challenging to capture the
variability of the process in every local structure in order to make
valid nonparametric inferences.




To circumvent the above difficulties, here we shall adopt a
different wild bootstrap technique. A similar technique was proposed
in Wu and Zhao (\citeyear{WZ07}) for nonparametric inference of the mean
function under stationary errors. The key idea is still uniform
Gaussian approximation of $\hat{Q}_{\alpha}(t)-Q_{\alpha}(t)$ on
$(0,1)$. However, instead of resorting to asymptotic theory, we
shall directly simulate the finite sample $\mathcal{L}^2$ and
$\mathcal
{L}^{\infty}$ norms of the Gaussian processes. More precisely, let us
assume the bandwidth $b_n(\alpha)=O(n^{-1/5})$ and
$Q_{\alpha}(t)\in\mathcal{C}^3[0,1]$ for illustrative purposes. Then by
the proofs of Theorems \ref{thm::scb} and \ref{thm::L2test}
in Section \ref{sec::pr}, we can obtain after elementary
calculations that
%
%
\begin{equation}\label{eq::bootsscb}
\sup_{t\in\mathcal{T}^*_n} \biggl|\frac{f(t,Q(t))}{\sigma(t)}[\tilde
{Q}(t)-Q(t)]- \mathfrak{X}_n(t) \biggr|=O_\mathbb{P}(n^{-11/20}\log^2 n)
\end{equation}
and
%
%
\begin{equation}\label{eq::bootsl2}
\biggl|\int_{\mathcal{T}_n^*}[\tilde{Q}_{\alpha}(t)-Q_{\alpha
}(t)]^2\tilde
{\pi}(t) \,dt-\int_{\mathcal{T}_n^*}\mathfrak{X}_n(t)^2 \,dt
\biggr|=O_\mathbb{P}
(n^\beta\log^{5/2} n),
\end{equation}
where
\[
\mathfrak{X}_n(t)=\sum_{i=1}^nV_iK^*_{b_n}(t_i-t)/(nb_n)
\]
with $(V_i)$ i.i.d. standard normal,
$\tilde{\pi}(t)=f^2(t,Q_{\alpha}(t))/\sigma^2(t)$ and $\beta=-0.95$.
Note that here the bias-corrected estimator $\tilde{Q}(\cdot)$ is
used in (\ref{eq::bootsscb}) in order to avoid estimating the
unpleasant bias term $\mu_2b_n^2Q''_{\alpha}(t)/2$ of the SCB.

An important observation of (\ref{eq::bootsscb}) and
(\ref{eq::bootsl2}) is that $\mathfrak{X}_n(t)$ does not depend on
the observations $(X_i)$ and has a simple and explicit form.
Therefore, one can generate a large sample of i.i.d. copies of
$\mathfrak{X}_n(t)$ and use the distributions of the $\mathcal
{L}^{\infty}$ and $\mathcal{L}^{2}$ norms of the sample to approximate
the distributions of the corresponding norms of $\tilde{Q}(t)-Q(t)$
by virtue of (\ref{eq::bootsscb}) and (\ref{eq::bootsl2}). The
following are the detailed procedures.
\begin{enumerate}[(6S)]
\item[(1)] Choose bandwidth $b_n$ according to the procedures in
Section \ref{sec::bandwidth}.
\item[(2)] Obtain $\tilde{Q}_{\alpha}(t)$ by (\ref{eq::ph-de}) and
(\ref{eq::jackknife}).
\item[(3)] Obtain estimate $\hat{\sigma}(t)$, $\hat{f}(t,Q_{\alpha
}(t))$ from (\ref{eq::var}) and (\ref{eq::den}) below. Let $\hat{\pi
}(t)=\hat{f}^2(t,Q_{\alpha}(t))/\hat{\sigma}^2(t)$.
\item[(4S)] Generate i.i.d. standard normal random variables $V_i$,
$i=1,2,\ldots,n$. Calculate $\varpi_{n,S}={\sup_{0\le t\le
1}}|\mathfrak{X}_n(t)|$.
\item[(4I)] Generate i.i.d. standard normal random variables $V_i$,
$i=1,2,\ldots,n$. Calculate $\varpi_{n,I}=\int_0^1\mathfrak
{X}^2_n(t) \,dt$.
\item[(5)] Repeat (4) $B$ times and obtain the estimated quantile
$\hat
{q}_{1-\beta}$ of $\varpi_{n}$.
\item[(6S)] The $100(1-\beta)\%$ SCB of $Q_{\alpha}(t)$ can be
constructed as $\tilde{Q}_{\alpha}(t)\pm\hat{q}_{1-\beta}/\sqrt
{\hat{\pi
}(t)}$. Accept the null hypothesis $Q_{\alpha}(\cdot)=Q^o_{\alpha
}(\cdot
)$ at level $1-\beta$ if and only $Q^o_{\alpha}(\cdot)$ is fully
contained in the SCB.
\item[(6I)] Accept the null hypothesis $Q_{\alpha}(\cdot
)=Q^o_{\alpha
}(\cdot)$ at level $1-\beta$ if and only if $\int_0^1[\tilde
{Q}_{\alpha
}(t)-Q^o_{\alpha}(t)]^2\hat{\pi}(t) \,dt\le\hat{q}_{1-\beta}$.
\end{enumerate}
One should use steps (1)--(3), (4S), (5) and (6S) when performing
hypothesis testing via the SCB. On the other hand, steps (1)--(3),
(4I), (5) and (6I) should be adopted when testing via the ISDT. The
number of replications $B$ can be chosen as 2000. It is immediate to
obtain $p$-values of the tests. For instance, the $p$-value of the
squared difference test is
$\mathbb{P}(\varpi_{n,I}>\int_0^1[\tilde{Q}_{\alpha
}(t)-Q^o_{\alpha
}(t)]^2\hat
{\pi}(t) \,dt)$,
which can be estimated by the bootstrap distribution of
$\varpi_{n,I}$.

Let $\{f(t,\theta)\}$ be a parametric family of functions that
depends on $t\in[0,1]$ and $\theta\in\Theta\subset\mathbb{R}^k$.
Often one
wants to test $Q_{\alpha}(t)=f(t,\theta)$ at level $1-\beta$ for
some unknown $\theta\in\Theta$. Under the null hypothesis, we have a
parametric model and one generally expects to obtain a root-$n$
consistent estimator $\hat{\theta}$ of the true parameter value
$\theta_0$ by the parametric quantile regression method of Keonker
(\citeyear{K05}). Note that the convergence rates of our SCB and ISDT tests
are always slower than $\sqrt{n}$. Therefore, if the null is true,
$f(t,\hat{\theta})$ can be treated as the true value of
$Q_{\alpha}(t)$ and one just needs to replace $Q^o_{\alpha}(\cdot)$
in steps (1)--(6) by $f(t,\hat{\theta})$ and the resulting testing
procedures are still valid.

The following proposition validates the above procedure for the case
$\{f(t,\theta)\}=\{\theta^{\top}\mathbf{g}(t)\}$, where $\theta\in
\mathbb{R}^k$
and $\mathbf{g}(t)\dvtx[0,1]\rightarrow\mathbb{R}^k$ is a known
function. We
shall first make the following constraint on $\mathbf{g}$:
\begin{enumerate}[(B1)]
\item[(B1)]
assume $\mathbf{g}(\cdot)\in\mathcal{C}[0,1]$ and $\mathcal
{G}:=\int
_0^1\mathbf{g}(t)\mathbf{g}^{\top}(t) \,dt$ is nonsingular.
\end{enumerate}
\begin{proposition}\label{prop::para_quantile}
Assume that $Q_{\alpha}(t)=\theta_0^{\top}\mathbf{g}(t)$ for some
$\theta_0\in\mathbb{R}^k$ and that conditions \textup{(A1)}, \textup{(A5)}
and \textup{(B1)} hold. Then
%
%
\begin{equation}\label{eq::prop21}
|\hat{\theta}_{\alpha}-\theta_0|=O_\mathbb{P}(n^{-1/2}),
\end{equation}
where
$\hat{\theta}_{\alpha}=\arg\min_{\theta}
\sum_{i=1}^n\rho_{\alpha}(X_i-\theta^{\top}\mathbf{g}(i/n))$.
\end{proposition}


\subsection{Estimation of the density and long-run variance
functions}\label{sec::estimate}
We see from (\ref{eq::bootsscb}) and (\ref{eq::bootsl2}) that
obtaining a good estimate of $f(t,Q(t))/\sigma(t)$ is necessary and
important for making our inferences. Here, we suggest using the
estimation techniques in Zhou and Wu (\citeyear{ZW09a}), which are essentially
local versions of the popular subsampling long-run variance
estimator and kernel density estimator for stationary data. Since
the time series is approximately stationary within comparatively
small time spans, the methods are shown to be consistent. See
Section 3.4 in Zhou and Wu (\citeyear{ZW09a}) for more details. For
the sake of
completeness, we present the estimators here.

For $t \in(0, 1)$, let $s_n(t) = \max(\lfloor n t - n b_n \rfloor,
1)$, $l_n(t) = \min(\lfloor n t + n b_n \rfloor, n)$ and
%
%
\begin{equation}\label{eq::Nnt}
\mathcal{N}_n(t)=\{i\in\mathbb{N}\dvtx s_n(t) \le i \le l_n(t)\}.
\end{equation}
Let $Z_{i,\alpha} = \psi_{\alpha}(X_i-\hat{Q}_{\alpha}(i/n))$. For a
sequence $m_n$ with $m_n \to\infty$ and $n b_n/m_n \to\infty$, we
shall estimate $\sigma^2(t)$ by
%
%
\begin{equation}\label{eq::var}
\hat{\sigma}^2(t)=\frac{m_n}{|\mathcal{N
}_n(t)|-m_n+1}\sum_{j=s_n(t)}^{l_n(t)-m_n+1}
\biggl(\frac{\sum_{i=j}^{j+m_n-1}Z_{i,\alpha}}{m_n}-\bar{Z}_n(t)
\biggr)^2,
\end{equation}
where $\bar{Z}_n(t) = \sum_{i \in\mathcal{N }_n(t) }Z_{i,\alpha} /
|\mathcal{N }_n(t)|$ and $|\mathcal{N}_n(t)| = l_n(t) - s_n(t) + 1$
is the
cardinality of $\mathcal{N}_n(t)$.

For $f(t,Q_{\alpha}(t))$, we shall use
%
%
\begin{equation}\label{eq::den}
\hat{f}(t,Q_{\alpha}(t))=\frac{1}{|\mathcal{N}_n(t)|h_n}
\sum_{i\in\mathcal{N}_n(t)}K^{\#}_{h_n}\bigl(\hat{Q}_{\alpha}(t)-X_i\bigr),
\end{equation}
where $K^{\#}\in\mathcal{K}$ is a kernel and $h_n$ is the bandwidth
satisfying $h_n \to0$ and $nb_nh_n \to\infty$.

We refer to Zhou and Wu (\citeyear{ZW09a}) for a discussion on the
selection of
the tuning parameters $m_n$ and $\tau_n$.

\subsection{Bandwidth selection}\label{sec::bandwidth}

Choosing a good bandwidth $b_n(\alpha)$ is important in practical
applications. For quantile curve estimation, Zhou and Wu (\citeyear{ZW09a})
proposed the following way to choose the bandwidth based on
modifications of existing bandwidth selectors for independent data.
By Theorem 1 of the latter paper regarding asymptotic normality of
$\hat{Q}_{\alpha}(t)$, we have
%
%
\begin{equation}\label{eq::band_ratio}
\frac{b^*_n(\alpha)}{b^{\mathrm{ind}}_n(\alpha)} =
\biggl[\frac{\int_0^1\sigma^2(t) \,dt}{\alpha(1-\alpha)} \biggr]^{1/5}
:=\rho^*(\alpha),
\end{equation}
where $b^*_n(\alpha)$ denotes the optimal weighted asymptotic mean
integrated\break squared error (AMISE) bandwidth, $b^{\mathrm{ind}}_n(\alpha)$ is
the optimal bandwidth obtained under independence and
$\rho^*(\alpha)$ is called the variance correction factor which
accounts for the dependence. Note that $\alpha(1-\alpha) =
\operatorname{Var}(J(t, Q_{\alpha}(t), \mathcal{F}_i))$. For independent
data, there
have been many discussions on bandwidth selection for nonparametric
quantile estimation; see, for instance, Yu and Jones (\citeyear
{YJ98}), Fan and
Gijbels (\citeyear{FG96}) and Ghosh and Draghicescu (\citeyear{GD02})
among others. Hence,
one could first select a bandwidth $b^{\mathrm{ind}}_n(\alpha)$ by treating
the data as if they were independent. After that, the variance
correction factor $\rho^*(\alpha)$ can be estimated by the following:
\[
\hat{\rho}^*(\alpha) = { { (\tilde{\sigma}^2)^{1/5}}
\over{ {(\alpha(1-\alpha))^{1/5} } } },
\]
where
\[
\tilde{\sigma}^2 = \frac{\tilde{m}}
{n-\tilde{m}+1}\sum_{j=1}^{n-\tilde{m}+1}
\Biggl(\frac{1}{\tilde{m}}
\sum_{i=j}^{j+\tilde{m}-1}
\varsigma_{i,\alpha}-\bar{\varsigma}_{n,\alpha} \Biggr)^2,
\]
$\varsigma_{i,\alpha} = \psi_{\alpha}(X_i -
\hat{Q}_{\alpha,b^{\mathrm{ind}}_n(\alpha)}(t))$, $\bar{\varsigma
}_{n,\alpha}
= \sum_{i=1}^n \varsigma_{i,\alpha} /n$ and $\tilde{m}=\lfloor n^{1/3}
\rfloor $. It can be shown that $\hat{\rho}^*(\alpha)$ is a consistent
estimate of $\rho^*(\alpha)$ and we shall point out that the
selected bandwidth $b^*_n(\alpha)$ typically varies with respect to
$\alpha$ [see also Yu and Jones (\citeyear{YJ98})]. Moreover, since the
jackknife bias reduction technique reduces the bias of the local
linear quantile estimates to second order, following Wu and Zhao
(\citeyear{WZ07}), we suggest using $b^{\mathrm{jack}}_n(\alpha)=2
b^*_n(\alpha)$ for
the nonparametric estimation when the jackknife is implemented. We
refer the interested reader to Sections 3.1.1 and 3.3 of Zhou and
Wu (\citeyear{ZW09a}) for more details on the bandwidth selection.

The bandwidth selected for quantile curve estimation provides a
reasonable starting point for nonparametric tests [Fan and Jiang
(\citeyear{FJ07})]. For the SCB test, we suggest using the bandwidth
$b^{\mathrm{jack}}_n(\alpha)$ following Eubank and Speckman (\citeyear{ES93})
and Wu and
Zhao (\citeyear{WZ07}). For the ISDT test $T_n^*$, Corollary
\ref{cor::optimality} implies that the bandwidth which renders the
optimal power is of order $n^{-2/9}=n^{-1/5}n^{-1/45}$. Following
Fan and Jiang (\citeyear{FJ07}), we suggest using the bandwidth
$b^{\mathrm{jack}}_n(\alpha)\times n^{-1/45}$ for the ISDT test.


\section{Simulation studies}\label{sec::simu}
In this section, we perform simulation studies to investigate the
accuracy and power of the proposed tests for moderate sample sizes.
Let us consider the following time-varying AR(1) model
%
%
\begin{equation}\label{eq::simu_model}
G(t,\mathcal{F}_i)=a_0(t)+a_1(t)G(t,\mathcal{F}_{i-1})+\delta
(t)\varepsilon_i,
\end{equation}
where $a_0(t)$, $a_1(t)$ and $\delta(t)$ are continuous functions on
$[0,1]$, ${\max_{t}}|a_1(t)|<1$, $\min_t\delta(t)>0$ and
$\varepsilon_i$ are i.i.d. with $\|\varepsilon_i\|_p<\infty$ for some
$p>0$. We observe the time series $X_i=G(i/n,\mathcal{F}_i)$,
$i=1,2,\ldots,n$.

It is clear that $G(t,\mathcal{F}_i)$ has the representation
%
%
\begin{equation}\label{eq::ar1}
G(t,\mathcal{F}_i)=\frac{a_0(t)}{1-a_1(t)}+\delta(t)\sum
_{j=0}^{\infty
}[a_1(t)]^{j}\varepsilon_{i-j}.
\end{equation}
The UGMC and local stationarity conditions of (\ref{eq::simu_model})
can be easily verified by the results in Section 4.1 of Zhou and Wu
(\citeyear{ZW09a}).

\subsection{Accuracy of the SCB and ISDT tests}
In this subsection, we describe a simulation study to compare the
accuracy of the asymptotic and bootstrap tests for both light tailed
and heavy tailed processes. To this end, we shall use the
time-varying AR(1) process (\ref{eq::simu_model}) with $a_0(t)=0$,
$a_1(t)=\sin(2\pi t)/2$ and $\delta(t)=\exp((t-1/4)^2)$. Consider
the following two scenarios:
\[
\mbox{(a) }\varepsilon_i\sim N(0,1);\qquad\mbox{(b)
}\varepsilon_i\sim S\alpha S(1.8).
\]
Here, $S\alpha S(1.8)$ denotes the standard symmetric $\alpha$ stable
distribution with index $1.8$ which has the characteristic function
$\exp(-|t|^{1.8})$, $t\in\mathbb{R}$. It is easy to show that the
$S\alpha
S(1.8)$ distribution has mean $0$ and infinite variance. Therefore,
scenarios (a) and (b) represent light tailed and heavy tailed
processes, respectively. Elementary calculations show that under
scenarios (a) and (b),
%
%
\begin{equation}\label{eq::quantile_simu}
Q_{\alpha}(t)=\frac{\delta(t)}{[1-|a_1(t)|^{\nu}]^{1/\nu
}}Q_{\alpha
}^{\varepsilon},
\end{equation}
where $Q_{\alpha}^{\varepsilon}$ is the $\alpha$th quantile of
$\varepsilon_i$ and $\nu=2$ and $1.8$ in scenarios (a) and (b),
respectively.

We shall compare type I error rates of the following six tests: the
asymptotic SCB test (AS) based on Theorem \ref{thm::scb}; the
asymptotic ISDT test (AI) based on Theorem \ref{thm::L2test}; the
bootstrap SCB test (BS); the bootstrap ISDT test (BI); the
asymptotic point-wise confidence band (PC) and the Bonferroni
confidence band (BF) based on Theorem 1 of Zhou and Wu (\citeyear
{ZW09a}). The
Bonferroni confidence band is simply the point-wise confidence band
at level $\beta/n$, where $\beta$ is the desired level. We generate
time-varying AR(1) processes under scenarios (a) and (b) with
length $n=300$ and perform the above six tests at the nominal level
$5\%$ for the following four quantile curves $\alpha=0.5$, $0.75$,
$0.9$ and $0.95$. Bandwidths are chosen according to Section
\ref{sec::bandwidth} and the critical values $\hat{q}_{0.95}$ of the
bootstrap tests are based on 2000 bootstrap samples. The simulated
type I error rates with 1000 replicates are shown in Table \ref{table1} below.

\begin{table}[b]
\tabcolsep=0pt
\caption{Simulated type \textup{I} error rates for the six tests with nominal
level $5\%$ under scenarios \textup{(a)} and \textup{(b)}.
Series length $n=300$ with $1000$ replicates}
\label{table1}
\begin{tabular*}{\tablewidth}{@{\extracolsep{\fill}}lcccccccc@{}}
\hline
& \multicolumn{2}{c}{$\bolds{\alpha=0.5}$} & \multicolumn{2}{c}{$\bolds{\alpha=0.75}$}
& \multicolumn{2}{c}{$\bolds{\alpha=0.9}$} & \multicolumn{2}{c@{}}{$\bolds{\alpha=0.95}$}\\
[-4pt]
& \multicolumn{2}{c}{\hrulefill} & \multicolumn{2}{c}{\hrulefill} & \multicolumn{2}{c}{\hrulefill}
& \multicolumn{2}{c@{}}{\hrulefill}\\
\textbf{Test} & \multicolumn{1}{c}{\textbf{(a)}} & \multicolumn{1}{c}{\textbf{(b)}}
& \multicolumn{1}{c}{\textbf{(a)}} & \multicolumn{1}{c}{\textbf{(b)}}
& \multicolumn{1}{c}{\textbf{(a)}} & \multicolumn{1}{c}{\textbf{(b)}}
& \multicolumn{1}{c}{\textbf{(a)}} & \multicolumn{1}{c@{}}{\textbf{(b)}} \\
\hline
AS & \phantom{0}2.2\%              & \phantom{0}2.6\%
& \phantom{0}2.7\%  & \phantom{0}2.8\%  & \phantom{0}2.1\%           & \phantom{0}9.9\%  & \phantom{0}2\%\phantom{0.}    & 21.7\% \\
AI & \phantom{0}8\%\phantom{0.}    & \phantom{0}6.5\%
& \phantom{0}7.2\%  & \phantom{0}6.9\%  & \phantom{0}7.4\%           & 10.4\%            & \phantom{0}6.7\%  & 21.4\%\\
BS & \phantom{0}4.8\%              & \phantom{0}6\%\phantom{0.}
& \phantom{0}4.4\%  & \phantom{0}5.4\%  & \phantom{0}4\%\phantom{0.} & 10.4\%            & \phantom{0}4\%\phantom{0.}    & 23.8\%\\
BI & \phantom{0}5.6\%              & \phantom{0}5.3\%
& \phantom{0}5.2\%  & \phantom{0}5.2\%  & \phantom{0}5.2\%           & \phantom{0}9.3\%  & \phantom{0}5.5\%  & 20.2\%\\
BF & \phantom{0}1.5\%              & \phantom{0}0.8\%
& \phantom{0}1.7\%  & \phantom{0}1.6\%  & \phantom{0}1.3\%           & \phantom{0}5.5\%  & \phantom{0}0.2\%  & 14.1\% \\
PC & 61.3\%                        & 61\%\phantom{0.}
& 53.9\%            & 53.1\%            & 45.4\%                     & 40.8\%            & 45.7\% & 44.2\%\\
\hline
\end{tabular*}
\end{table}

It is clear from the output of Table \ref{table1} that point-wise confidence
bands are not appropriate for nonparametric inference. The
Bonferroni confidence band test is too conservative, namely the band
is too wide. On the other hand, the asymptotic SCB test (AS) is
conservative and the asymptotic ISDT test (AI) tends to slightly
inflate the type I error. As discussed in Section \ref{sec::boots},
the loss of accuracy of the asymptotic tests is due to their slow
convergence rates.

For the bootstrap tests, the nominal type I error is achieved except
for extreme quantiles of heavy tailed processes. It is not difficult
to see that data is relatively sparse at high quantiles and very
large jumps occur more frequently for heavy tailed processes. These
facts suggest that for extreme quantile inference of heavy tailed
processes, a relatively large sample size is needed in order to
achieve the desired accuracy. We performed the bootstrap tests at
$5\%$ level for scenario (b) for $n=500$ and with 1000 replicates.
Simulated type I errors of the BS and BI tests were $0.061$ and
$0.056$, respectively, for the $90\%$ quantile curve. However, for the
$95\%$ quantile, the corresponding simulated type I errors were
$0.139$ and $0.131$, respectively, which were still way larger than
the nominal. When $n$ was increased to $1000$, accuracy of the
bootstrap tests were achieved for the $95\%$ quantile curve under
scenario (b).

\subsection{Power comparison of the SCB and ISDT tests}\label{sec::simu_power}
\begin{figure}

\includegraphics{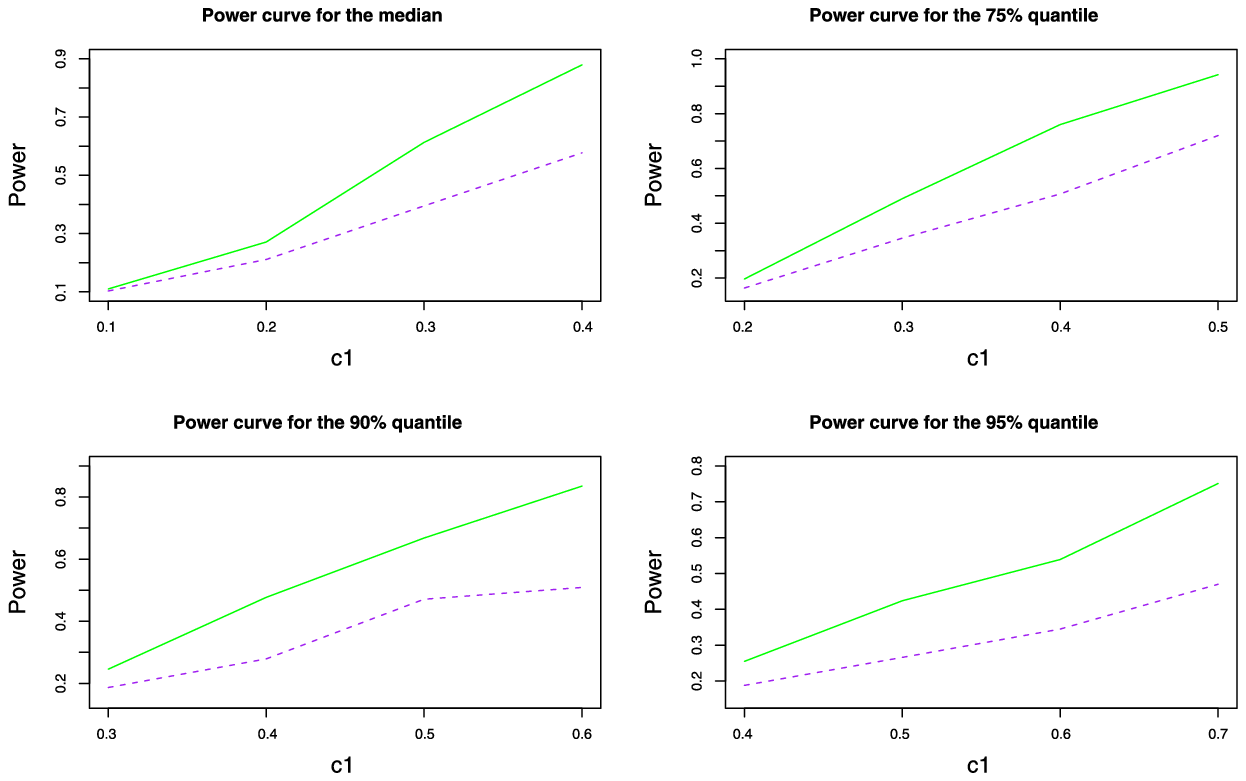}

\caption{Power curves for the $50\%$, $75\%$,
$90\%$ and $95\%$ quantiles under case \textup{(i)} of Section
\protect\ref{sec::simu_power}. The solid lines are the power curves
for the
ISDT test and the dashed lines are the power curves for the SCB
test.}\label{fid:powerflat}\vspace*{20pt}
%

\includegraphics{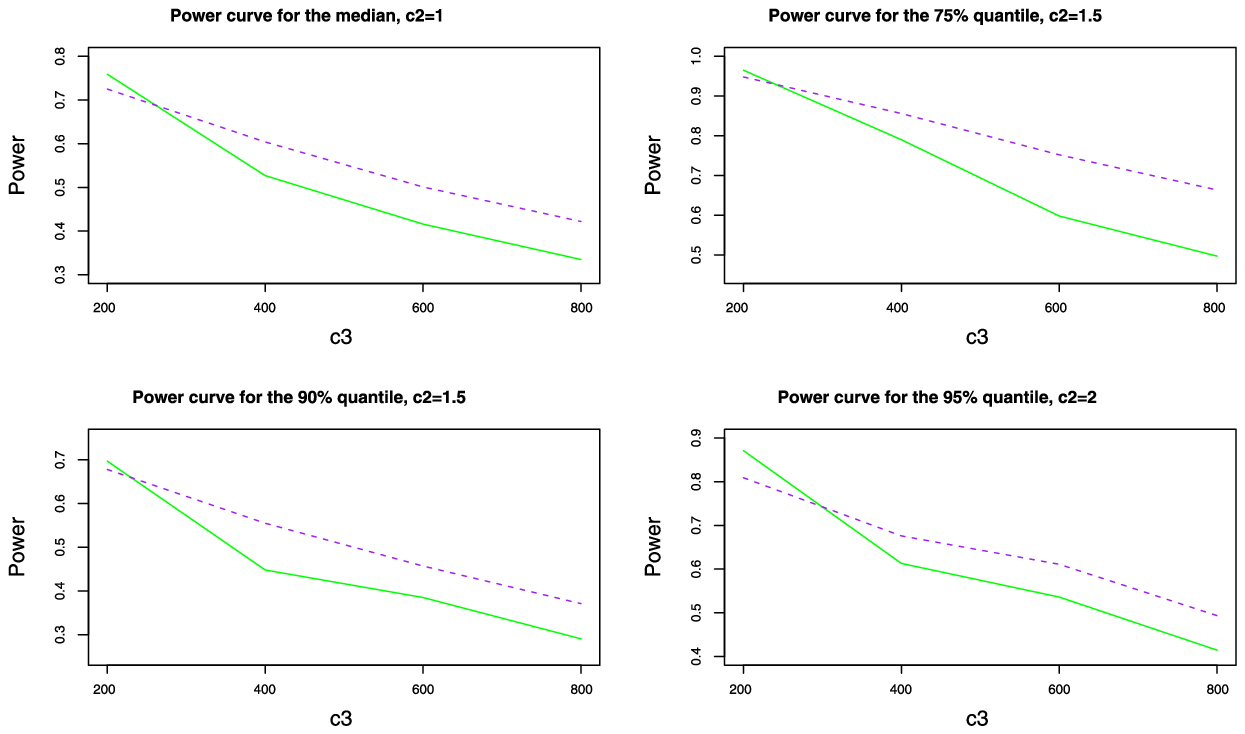}

\caption{Power curves for the $50\%$, $75\%$,
$90\%$ and $95\%$ quantiles under case \textup{(ii)} of Section
\protect\ref{sec::simu_power}. The solid lines are the power curves
for the
ISDT test and the dashed lines are the power curves for the SCB
test.}\label{fid:powerbumpy}
\end{figure}

As discussed in Remark~\ref{rem::1}, for a moderate sample size,
power of the SCB and ISDT tests are greatly influenced by the shape
of $Q_{\alpha}(t)-Q_{\alpha}^o(t)$. Recall that $Q_{\alpha}^o(t)$ is
the hypothesized value of $Q_{\alpha}(t)$. In this subsection, we
present a simulation to compare the power performance of the above
two tests under various shapes of $Q_{\alpha}(t)-Q_{\alpha}^o(t)$.
To this end, we consider model (\ref{eq::simu_model}) with
$a_0(t)=\varphi(t)(1-a_1(t))$, $a_1(t)=\sin(2\pi t)/2$,
$\delta(t)=\exp((t-1/4)^2)$ and $\varepsilon_i$ i.i.d. $N(0,1)$. Then
%
%
\begin{equation}\label{eq::quant_simu2}
Q_{\alpha}(t)=\varphi(t)+Q_{\alpha}^o(t)\qquad \mbox{where } Q_{\alpha
}^o(t)=\frac{\delta(t)}{[1-a^2_1(t)]^{1/2}}Q_{\alpha}^{\varepsilon}.
\end{equation}
We test the hypothesis $H_0\dvtx
Q_{\alpha}(\cdot)=Q_{\alpha}^o(\cdot)$ versus $H_a\dvtx
Q_{\alpha}(\cdot)\neq Q_{\alpha}^o(\cdot)$. Consider the following
two situations:
\[
\mbox{(i) } \varphi(t)=c_1;\qquad \mbox{(ii) }\varphi(t)=c_2 \exp
\bigl(-c_3(t-1/2)^2\bigr),
\]
where $c_i$, $i=1,2,3$, are positive constants. Cases (i) and
(ii) correspond to flat and bumpy differences of the true and
hypothesized quantile curves, respectively. Note that as $c_3$ gets
larger, we observe sharper peaks in $\varphi(t)$. In our simulations,
we follow steps (1) to (6) in Section \ref{sec::boots} and
perform the SCB and ISDT tests at $5\%$ level with $n=300$ and
$\alpha=0.5$, $0.75$, $0.9$ and $0.95$. Various values of $c_1$ and
$(c_2,c_3)$ are investigated and for each of the values we perform
1000 replicates and record the simulated probability of rejecting
the null hypothesis. The simulated power curves of case (i) and
(ii) are shown in Figures \ref{fid:powerflat} and
\ref{fid:powerbumpy}, respectively.

%

Generally, the displays in Figures \ref{fid:powerflat} and
\ref{fid:powerbumpy} are consistent with our arguments in Remark
\ref{rem::1} that for a moderate sample size the ISDT test is more
powerful when $\varphi(\cdot)$ is flat and the SCB test performs
better when $\varphi(\cdot)$ changes abruptly. On the other hand,
when the peak of $\varphi(\cdot)$ is not sharp enough; namely when
$c_3$ is relatively small, we see from Figure \ref{fid:powerbumpy}
that the ISDT test is still more powerful than the SCB test, which
can be explained by the fact that the ISDT test asymptotically
dominates the SCB test.

\section{The global tropical cyclone data}\label{sec::data}
\begin{figure}[b]

\includegraphics{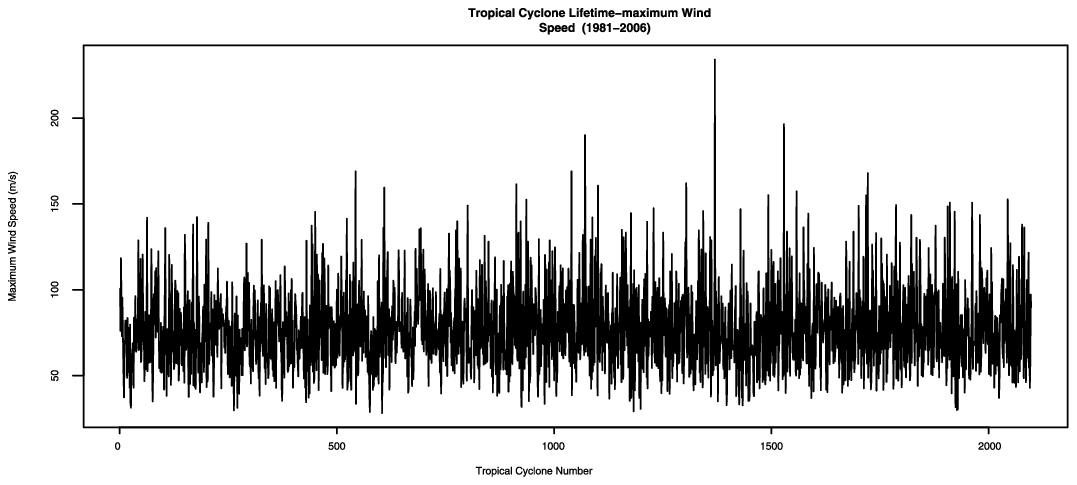}

\caption{Time series plot of tropical cyclone lifetime-maximum wind
speed (1981--2006).}\label{fid:wind}
\end{figure}

One of the most important consequences of global warming is the
increase of ocean temperature. Theoretical arguments and modeling
studies indicate that tropical cyclone winds should increase with
increasing ocean temperature [Elsner, Kossin and Jagger (\citeyear
{EKJ08})]. Meanwhile,
climatologists are very interested in finding empirical evidences on
the change of intensity of tropical cyclone winds. In a recent
paper, Elsner, Kossin and Jagger (\citeyear{EKJ08}) tackled the latter
problem partially by
fitting \textit{linear} trends for quantiles of the global tropical
cyclone data. The data set contains satellite-derived
lifetime-maximum wind speeds of 2098 tropical cyclones over the
globe during the period 1981--2006. It is available at James Elsner's
website at
\href{http://myweb.fsu.edu/jelsner/extspace/globalTCmax4.txt}{http://myweb.fsu.edu/jelsner/extspace/}
\href{http://myweb.fsu.edu/jelsner/extspace/globalTCmax4.txt}{globalTCmax4.txt}.
We shall refer to Elsner, Kossin and Jagger (\citeyear{EKJ08}) for a
detailed description on
how the data are obtained and the related issues. Figure
\ref{fid:wind} shows the time series plot of the data. Significant
linearly increasing trends were found in high quantiles of the
global tropical cyclone data in Elsner, Kossin and Jagger (\citeyear
{EKJ08}). In other
words, the worst tropical cyclones are getting stronger over the
globe.

%

In this section, we are mainly interested in the following issues.
First, we shall compare our quantile-based tests with the mean-based
approaches to see whether the increase of intensity of the worst
tropical cyclones is due to a shift of mean. Second, we shall test
whether the linear trend used in Elsner, Kossin and Jagger (\citeyear
{EKJ08}) is sufficient
to describe the change of the high quantiles of the global cyclone
winds. This is an important issue because the parametric linearity
assumption implies homogeneous change of the wind speeds. However,
much complex yet important information on the dynamics of high wind
speeds may be buried under this assumption.

%
\begin{figure}

\includegraphics{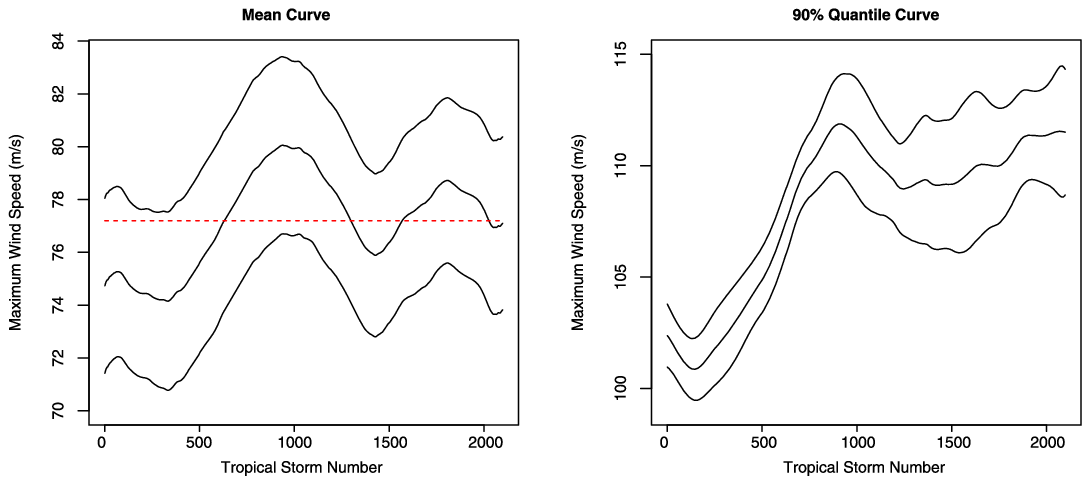}

\caption{Trends in mean and $90\%$ quantile of global tropical cyclone
lifetime-maximum wind
speed, with $95\%$ simultaneous confidence bands.}\label{fid:ccc}
\end{figure}

To address the first issue, we follow the procedures in Zhou and Wu
(\citeyear{ZW09b}) and provide a SCB for the mean curve of the global tropical
cyclone series. The bandwidth is chosen as 0.16. The left panel of
Figure \ref{fid:ccc} shows the 95\% SCB. The height of the
horizontal line is the average of all cyclone wind speeds. It is
clear that the horizontal line is fully contained in the SCB. The
$p$-value for testing constancy of the mean trend is 0.11. Therefore,
statistically there is no trend in the average tropical cyclone
speeds. In fact, one of the main reasons why scientists did not find
interesting signals in the intensity of global tropical cyclone
winds before Elsner, Kossin and Jagger (\citeyear{EKJ08}) is due to
the focus on mean
trends.

We performed the SCB and ISDT tests on the linearity of the 90\%
quantile curve to address the second issue. The bandwidths of the
SCB and ISDT tests are chosen as 0.21 and 0.18, respectively. Both
tests give $p$-values less than 0.001. The right panel of Figure
\ref{fid:ccc} shows the 95\% SCB. It is clear from the SCB that
there is an inhomogeneous increasing trend of the 90\% quantile
curve. More specifically, the quantile curve underwent a sharp
increase during the period November 1982 to September 1992, and
after that the trend became flat. This information cannot be
retrieved if the linear trend analysis in Elsner, Kossin and Jagger
(\citeyear{EKJ08}) is
adopted. Furthermore, the inhomogeneous trend in high quantiles
provides climatologists with useful information on the underlying
complex mechanisms producing the hurricane climate.


\section{Proofs}
\label{sec::pr} Unless otherwise specified, we will only prove
results for $\alpha=1/2$, since results for other quantiles follow
similarly. We shall also omit the subscript $\alpha$ in the notation
if no confusion will be caused.

Consider a system $H(t,x,\mathcal{F}_i)$ defined in Section
\ref{sec::model}. Let $g\dvtx[0,1]\mapsto\mathbb{R}$ be a measurable
function. Define $\mathfrak{H}(t,\mathcal{F}_i)=H(t,g(t),\mathcal
{F}_i)$. Then
$Z_i:=\mathfrak{H}(t_i,\mathcal{F}_i)$ $i=1,2,\ldots,n$, defines a
nonstationary time series. Recall $t_i=i/n$. Let
$S_Z(i)=\sum_{j=1}^iZ_j$, $i=1,2,\ldots,n$. The following invariance
principle for $(Z_i)_{i=1}^n$ plays an important role in both
establishing the asymptotic theory of the testing methods and
justifying the wild bootstrap procedures.
\begin{theorem}\label{thm::sip}
Assume that \textup{(i)} $\mathbb{E}\mathfrak{H}(t,\mathcal{F}_i)=0$ for all
$t\in[0,1]$;
that \textup{(ii)} ${\sup_{(t,x)\in[0,1]\times\mathbb{R}}}\|H(t,x,\mathcal
{F}_i)\|_4<c_0$ for
some finite constant $c_0$; that \textup{(iii)} $H$ satisfies $\operatorname{UGMC}(4)$ and
that \textup{(iv)} there exists $q\ge1/4$, such that for all $i\ge0$
%
%
\begin{equation}\label{eq::lipsip}
\| \mathcal{P}_0\{\mathfrak{H}(s_1,\mathcal{F}_i) - \mathfrak
{H}(s_2,\mathcal{F}
_i) \}
\| \le C
|s_1-s_2|^q
\end{equation}
holds for all $s_1, s_2 \in[0,1]$, where $\mathcal{P}_k(\cdot
)=\mathbb{E}
(\cdot
|\mathcal{F}_k)-\mathbb{E}(\cdot|\mathcal{F}_{k-1})$ for $k\ge0$. Then
on a richer probability space, there exists i.i.d. standard normal
random variables $(V_i)_1^n$ and a process $S_Z^\circ(i)$ with
$\{S^\circ_Z(i) \}_{i=1}^n \stackrel{\mathcal{D}}{=}\{S_Z(i) \}
_{i=1}^n$, such that
%
%
\begin{equation}\label{eq::sip}
\max_{i\le n}\Biggl|S^\circ_Z(i)-\sum_{j=1}^i\upsilon(j/n)V_j\Biggr|=o_\mathbb{P}
(n^{1/4}\log^2 n),
\end{equation}
where
$\upsilon(s)= [\sum_{k\in\mathbb{Z}}\operatorname{cov}(H(s,g(s),\mathcal
{F}_0),H(s,g(s),\mathcal{F}_k)) ]^{1/2}, $ $0\le s\le1$.
\end{theorem}
\begin{pf}
Note $\upsilon(s)=\|\mathcal{P}_0\sum_{k=0}^{\infty}\mathfrak
{H}(s,\mathcal
{F}_k)\|$. Let $\mathfrak{H}(t,\mathcal{F}_i)=0$ if $t>1$. By Corollary
1 of Wu and Zhou (\citeyear{WZ09}), it follows that under conditions
(i) and
(ii), there exist i.i.d. standard normal random variables $(V_i)_1^n$
and a process $S_Z^\circ(i)$ with $\{S^\circ_Z(i)\}_{i=1}^n \stackrel
{\mathcal{D}}{=}
\{S_Z(i) \}_{i=1}^n$, such that
%
%
\begin{equation}\label{eq::thm11}
\max_{i\le n}\Biggl|S^\circ_Z(i)-\sum_{j=1}^i\tilde{\upsilon
}_jV_j\Biggr|=o_\mathbb{P}
(n^{1/4}\log^{3/2} n),
\end{equation}
where $\tilde{\upsilon}_{j}=\|\mathcal{P}_0\sum_{k=0}^{\infty
}\mathfrak
{H}(t_j+t_k,\mathcal{F}_k)\|$. Based on (\ref{eq::lipsip}),
%
%
\begin{equation}\label{eq::thm13}
\|\mathcal{P}_0\{\mathfrak{H}(t_i,\mathcal{F}_k)-\mathfrak
{H}(t_i+t_k,\mathcal
{F}_k)\}\|\le t^{q}_k
\end{equation}
for $0\le k\le n-i$. On the other hand, by Theorem 1 of Wu (\citeyear
{W05}), we obtain
%
%
\begin{equation}\label{eq::thm1a}
\|\mathcal{P}_0\{\mathfrak{H}(t_i,\mathcal{F}_k)-\mathfrak
{H}(t_i+t_k,\mathcal
{F}_k)\}\|\le2\delta_H(k,2)
\end{equation}
for all $k\ge0$. By (\ref{eq::thm13}), (\ref{eq::thm1a}) and the fact
that $\delta_H(k,2)\le\delta_H(k,4)=O(\chi^k)$ for some $0<\chi
<1$, we
have for all $1\le j\le n-\lceil n^{1/4}\rceil $,
\begin{eqnarray*}
\bigl(\upsilon(j/n)-\tilde{\upsilon}_{j}\bigr)^2&\le& \Biggl\|\mathcal{P}_0\sum
_{k=0}^{\infty
}[\mathfrak{H}(t_j,\mathcal{F}_k)-\mathfrak{H}(t_j+t_k,\mathcal
{F}_k)]\Biggr\|
^2\\
&\le& \Biggl(\sum_{k=0}^{\infty}\|\mathcal{P}_0[\mathfrak
{H}(t_j,\mathcal
{F}_k)-\mathfrak{H}(t_j+t_k,\mathcal{F}_k)]\|\Biggr)^2\\
&\le& \Biggl(\sum_{k=0}^{k^*}t^{q}_k+\sum_{k=k^*}^{\infty}2\chi
^k\Biggr)^2=O(n^{-1/2}\log n^{5/2}),
\end{eqnarray*}
where $k^*=\lfloor -\log n/\log\chi\rfloor $.
Note for $n-\lceil n^{1/4}\rceil <j\le n$, $(\upsilon(j/n)-\tilde
{\upsilon
}_{j})^2=O(1)$. Therefore,
$\sum_{j=1}^n(\upsilon(j/n)-\tilde{\upsilon}_{j})^2=O(n^{1/2}\log^{5/2}
n)$. Hence,
%
%
\begin{eqnarray}\label{eq::thm12}
\max_{i\le n}\Biggl|\sum_{j=1}^i\upsilon(j/n)V_j-\sum_{j=1}^i\tilde
{\upsilon
}_{j}V_j\Biggr|&=&O_\mathbb{P}((n^{1/2}\log^{5/2} n)^{1/2}\log^{1/2}
n)\nonumber\\[-8pt]\\[-8pt]
&=&o_\mathbb{P}(n^{1/4}\log^2 n).\nonumber
\end{eqnarray}
By (\ref{eq::thm11}) and (\ref{eq::thm12}), Theorem \ref{thm::sip} follows.
\end{pf}
\begin{remark}
The Gaussian approximation in Theorem \ref{thm::sip} shows that
partial sums of a short range dependent (SRD) locally stationary
process can be approximated by weighted sums of i.i.d. standard normal
random variables. We shall call
$\upsilon^2(t)=\sum_{k\in\mathbb{Z}}\operatorname{cov}(\mathfrak
{H}(t,\mathcal
{F}_0),\mathfrak{H}(t,\mathcal{F}_k))$ the long-run variance of the
system $\{\mathfrak{H}(t,\mathcal{F}_i)\}$ at time $t$. The weight
$\upsilon(t)$ captures the local dependence of the series $(Z_i)$ at
$t$; while the fluctuation of $\upsilon(t)$ on $[0,1]$ is due to the
nonstationarity of the series. Following the arguments in Wu and
Zhou (\citeyear{WZ09}), the bound $o_\mathbb{P}(n^{1/4}\log^2 n)$ in
(\ref{eq::sip}) is
optimal within a multiplicative logarithmic factor.
\end{remark}

We shall state and prove the following Lemma \ref{lem::approximation}
and Lemma \ref{lem::guassion_max} before we proceed to the proof of
Theorem \ref{thm::scb}.
\begin{lemma}\label{lem::approximation}
Assume $Q(t)\in\mathcal{C}^2[0,1]$; that conditions \textup{(A1)--(A6)} and
\textup{(K1)}
hold, and that $b_n\rightarrow0$ with $nb^{3/2}_n\rightarrow\infty$,
then we have
\begin{eqnarray*}
&&{\sup_{t\in\mathcal{T}_n}}|f(t,Q(t))[\hat{Q}(t)-Q(t)]-
S_n(t)-\mathbb{E}
S_n^*(t)|\\
&&\qquad=O_{\mathbb{P}} \biggl( {{\pi_n^{1/2}\log n+b_n^{3/4}} \over\sqrt{nb_n}}+
b_n\pi_n+ \pi_n^{2} \biggr),
\end{eqnarray*}
where $\pi_n=(nb_n)^{-1/2}(\log n+(b_n)^{-1/4}+(nb^5_n)^{1/2})$,
%
%
\begin{equation}\label{eq::lem11}
S_n(t)= \sum_{i=1}^n\psi\bigl(X_i-Q(t_i)\bigr)K_{b_n}(t-t_i)/(nb_n)
\end{equation}
with $\psi_{\alpha}(x)=\alpha-I\{x\le0\}$, and
\[
S_n^*(t)=\sum_{i=1}^n\psi\bigl(X_i-Q(t)-(t_i-t)Q'(t)\bigr)K_{b_n}(t-t_i)/(nb_n).
\]
\end{lemma}
\begin{pf}
A careful check of the proof of Theorem 3 of Zhou and Wu (\citeyear{ZW09a})
shows that under conditions of Lemma \ref{lem::approximation}
%
%
\begin{eqnarray}\label{eq::lem12}
&&{\sup_{t\in\mathcal{T}_n}}|f(t,Q(t))[\hat{Q}(t)-Q(t)]-
S_n^*(t)|\nonumber\\[-8pt]\\[-8pt]
&&\qquad=O_{\mathbb{P}} \biggl( {{\pi_n^{1/2}\log n} \over\sqrt{nb_n}}+
b_n\pi_n+ \pi_n^{2} \biggr).\nonumber
\end{eqnarray}
Let $\psi(i,t)=\psi(X_i-Q(t_i))-\psi(X_i-Q(t)-(t_i-t)Q'(t))$,
\begin{eqnarray*}
M_n(t) &=& \sum_{i=1}^n\mathcal{P}_i\psi(i,t)K_{b_n}(t-t_i)/(nb_n),
\\
N_n(t) &=& \sum_{i=1}^n\{\mathbb{E}[\psi(i,t)|\mathcal
{F}_{i-1}]-\mathbb{E}\psi(i,t)\}
K_{b_n}(t-t_i)/(nb_n).
\end{eqnarray*}
Note the summands of $M_n(t)$ form a triangular array of martingale
differences and $N_n(t)$ is differentiable with respect to $t$.
Using similar chaining arguments as those in the proof of Lemmas 5
and 6 in Zhou and Wu (\citeyear{ZW09a}), we have
%
%
\begin{equation}\label{eq::lem13}
{\sup_{t\in[0,1]}}|M_n(t)|=O_\mathbb{P}\biggl(\frac{b_n\log n}{\sqrt
{nb_n}}\biggr),\qquad
{\sup_{t\in[0,1]}}|N_n(t)|=O_\mathbb{P}\biggl(\frac{b_n^{3/4}}{\sqrt{nb_n}}\biggr).
\end{equation}
Note $S_n(t)-(S_n^*(t)-\mathbb{E}S_n^*(t))=M_n(t)+N_n(t)$. Therefore, by
(\ref{eq::lem12}), (\ref{eq::lem13}) and the fact that $b_n\log
n=o(b_n^{3/4})$, we show that Lemma \ref{lem::approximation} holds.
\end{pf}
\begin{lemma}\label{lem::guassion_max}
Let $\{V_i\}$ be i.i.d. standard normal random variables. Assume
condition \textup{(K1)} and $b_n\rightarrow0$ with $nb_n/(\log n)^2\rightarrow
\infty$. Then for any $x\in\mathbb{R}$
\[
\lim_{n\rightarrow\infty}\mathbb{P}\Biggl[\frac{1}{\sqrt{\phi
nb_n}}\sup_{t\in
\mathcal{T}_n}\Biggl|\sum_{i=1}^nV_iK_{b_n}(t_i-t)\Biggr|
-B(m^*)\le\frac{x}{\sqrt{2\log m^*}} \Biggr]= e^{-2e^{-x}}.
\]
\end{lemma}

Lemma \ref{lem::guassion_max} follows from classic results for extremes
of Gaussian processes. See, for example, Bickel and Rosenblatt
(\citeyear{BR73})
and Lemma 2 of Wu and Zhao (\citeyear{WZ07}). Details are omitted.
\begin{pf*}{Proof of Theorem \ref{thm::scb}}
Consider the system $\{1/2-J(t,x,\mathcal{F}_i)\}$, $(t,x)\in
[0,1]\times\mathbb{R}$.
Recall $J(t,x,\mathcal{F}_i)=I\{G(t,\mathcal{F}_i)\le x\}$. Let
$\mathfrak{J}(t,\mathcal{F}_i)=1/2-J(t$,\break $Q(t),\mathcal{F}_i)$. We
shall first show
that $\mathfrak{J}(t,\mathcal{F}_i)$ satisfies conditions of Theorem
\ref{thm::sip}. Obviously $\mathbb{E}\mathfrak{J}(t,\mathcal
{F}_i)=0$ and
$\|1/2-J(t,x,\mathcal{F}_i)\|_4\le1$ for all $t\in[0,1]$ and $x\in
\mathbb{R}$.
Based on (A5), condition (iii) of Theorem \ref{thm::sip} holds. We
now check (iv). Note for any $i\ge0$
%
%
\begin{eqnarray}\label{eq::thm14}\quad
\|\mathcal{P}_0\{\mathfrak{J}(s_1,\mathcal{F}_i) - \mathfrak
{J}(s_2,\mathcal{F}_i)
\}\|
&\le& \|\mathfrak{J}(s_1,\mathcal{F}_i) - \mathfrak{J}(s_2,\mathcal
{F}_i)\|
\nonumber\\[-8pt]\\[-8pt]
&=&\|J(s_1,Q(s_1),\mathcal{F}_i) - J(s_2,Q(s_2),\mathcal{F}_i)
\|.\nonumber
\end{eqnarray}
Based on (A6) and the smoothness condition on $Q(t)$,
%
%
\begin{equation}\label{eq::thm15}
\|J(s_2,Q(s_1),\mathcal{F}_i)-J(s_2,Q(s_2),\mathcal{F}_i)\|\le
C|s_1-s_2|^{1/2}.
\end{equation}
On the other hand,
%
%
\begin{equation}\label{eq::thm16}
\|J(s_1,Q(s_1),\mathcal{F}_i)-J(s_2,Q(s_1),\mathcal{F}_i)\|\le I+\mathit{I
I},
\end{equation}
where $I=\|[J(s_1,Q(s_1),\mathcal{F}_i)-J(s_2,Q(s_1),\mathcal
{F}_i)]I\{|\zeta
_i(s_1)-\zeta_i(s_2)|\le\delta\}\|$ and
$\mathit{I I}=\{|\zeta_i(s_1)-\zeta_i(s_2)|>\delta\}\|$. Recall $\zeta
_i(t)=G(t,\mathcal{F}_i)$. Using condition (A6), we have for all
$\delta>0$,
%
%
\begin{eqnarray}\label{eq::thm17}
I&\le&\|I\{Q(s_1)-\delta\le\zeta_i(s_1)\le Q(s_1)\}\|+\|I\{
Q(s_1)\le
\zeta_i(s_1)\le Q(s_1)+\delta\}\|\hspace*{-32pt}\nonumber\\[-8pt]\\[-8pt]
&\le& C\delta^{1/2}.\nonumber
\end{eqnarray}
By condition (A3), there exists $q\ge1$, such that
%
%
\begin{equation}\label{eq::thm18}
\mathit{I I}\le\frac{\|\zeta_i(s_1)-\zeta_i(s_2)\|_q^{q/2}}{\delta
^{q/2}}\le
C\frac{|s_1-s_2|^{q/2}}{\delta^{q/2}}.
\end{equation}
Let $\delta=|s_1-s_2|^{q/(q+1)}$, then (\ref{eq::thm16}), (\ref
{eq::thm17}) and (\ref{eq::thm18}) imply that
%
%
\begin{equation}\label{eq::thm19}
\|J(s_1,Q(s_1),\mathcal{F}_i)-J(s_2,Q(s_1),\mathcal{F}_i)\|\le C|s_1-s_2|^s,
\end{equation}
where $s=\frac{q}{2(1+q)}\ge\frac{1}{4}$. By (\ref{eq::thm14}),
(\ref{eq::thm15}) and (\ref{eq::thm19}), we have condition (iv) of
Theorem \ref{thm::sip} holds. Therefore, Theorem \ref{thm::sip}
implies that there exist i.i.d. standard normal random variables
$(V_i)_1^n$, such that
%
%
\begin{equation}\label{eq::thm110}
\max_{i\le n}\Biggl|S_{\mathfrak{J}}(i)-\sum_{j=1}^i\sigma
(t_j)V_j\Biggr|=o_\mathbb{P}
(n^{1/4}\log^2 n),
\end{equation}
where $S_{\mathfrak{J}}(i)=\sum_{k=1}^i\mathfrak{J}(t_k,\mathcal
{F}_k)$. Recall
$\sigma(t)$ was defined in (\ref{eq::lrv}). Define
%
%
\begin{eqnarray}\label{eq::thm111}
\Xi_n(t)&=&\sum_{i=1}^n\sigma(t_i)V_iK_{b_n}(t_i-t)/(nb_n),
\\
%
%
\label{eq::thm112}
\Omega_n(t)&=&\Biggl\{K_{b_n}(t_1-t)+\sum
_{i=2}^n|K_{b_n}(t_i-t)-K_{b_n}(t_{i-1}-t)|\Biggr\}\bigg/(nb_n).
\end{eqnarray}
By the summation by parts formula and (\ref{eq::thm110}), we obtain
%
%
\begin{eqnarray}\label{eq::thm113}
{\sup_{t\in\mathcal{T}_n}}|S_n(t)-\Xi_n(t)|&\le&\sup_{t\in\mathcal
{T}_n}\Omega_n(t)\max_{i\le n}\Biggl|S_{\mathfrak{J}}(i)-\sum
_{j=1}^i\sigma
(t_j)V_j\Biggr|\nonumber\\[-8pt]\\[-8pt]
&=&o_\mathbb{P}\biggl(\frac{n^{1/4}\log^2 n}{nb_n} \biggr).\nonumber
\end{eqnarray}
By the Lipschitz continuity of $\sigma(t)$ in (A2), it is easy to see that
%
%
\begin{equation}\label{eq::thm114}
{\sup_{t\in\mathcal{T}_n}}|\Xi_n(t)-\Xi^*_n(t)|=O_\mathbb{P}\biggl(\frac
{b_n\log
n}{\sqrt{nb_n}} \biggr),
\end{equation}
where $\Xi^*_n(t)=\sigma(t)\sum_{i=1}^nV_iK_{b_n}(t_i-t)/(nb_n)$.
By Lemma \ref{lem::approximation}, (\ref{eq::thm113}), (\ref
{eq::thm114}) and the fact that $b_n\log n=o(b_n^{3/4})$, it follows that
%
%
\begin{eqnarray}\label{eq::thm115}
&&{\sup_{t\in\mathcal{T}_n}}|f(t,Q(t))[\hat{Q}(t)-Q(t)]- \Xi
^*_n(t)-\mathbb{E}
S_n^*(t)|\nonumber\\[-8pt]\\[-8pt]
&&\qquad=O_\mathbb{P}\biggl( \frac{\pi_n^{1/2}\log n+b_n^{3/4}}{\sqrt
{nb_n}}+b_n\pi
_n+\pi_n^{2}+\frac{n^{1/4}\log^2 n}{{nb_n}} \biggr).\nonumber
\end{eqnarray}
It is easy to check that under bandwidth conditions of Theorem
\ref{thm::scb}, the right-hand side of (\ref{eq::thm115}) is of
order $o_\mathbb{P}((nb_n\log n)^{-1/2})$. On the other hand, by condition
(A1) and a Taylor expansion, we have
%
%
\begin{eqnarray}\label{eq::thm116}
{\sup_{t\in\mathcal{T}_n}}|\mathbb{E}S_n^*(t)-\mu
_2f(t,Q(t))Q''(t)b^2_n/2|&=&O \biggl(b_n^3+\frac{1}{nb_n}
\biggr)\nonumber\\[-8pt]\\[-8pt]
&=&o((nb_n\log n)^{-1/2}).\nonumber
\end{eqnarray}
By Lemma \ref{lem::guassion_max}, (\ref{eq::thm115}) and (\ref
{eq::thm116}), Theorem \ref{thm::scb} follows.
\end{pf*}

Let
%
%
\begin{eqnarray}\label{eq::lem33}
\mathcal{S}_n(t)&=& \sum_{i=1}^n\psi
\bigl(X_i-Q(t_i)\bigr)K^*_{b_n}(t-t_i)/(nb_n),\nonumber\\[-8pt]\\[-8pt]
\mathcal{S}_n^o(t)&=&\sum_{i=1}^n\psi
\bigl(X_i-Q(t)-(t_i-t)Q'(t)\bigr)K^*_{b_n}(t-t_i)/(nb_n).\nonumber
\end{eqnarray}
Define
%
%
\begin{equation}\label{eq::ynt}
\mathfrak{Y}_n(t)=\tilde{Q}(t)-Q(t)-\mathbb{E}\mathcal{S}_n^o(t)/f(t,Q(t)).
\end{equation}
We shall introduce several lemmas before proceeding to the proof of
Theorem \ref{thm::L2test}.
%
%
\begin{lemma}\label{lem::linearform}
Under the conditions of Theorem \ref{thm::L2test}, we have
%
%
\begin{equation}\label{eq::lem31}
\sqrt{n}\int_{\mathcal{T}_n^*}\mathfrak{Y}_n(t)\pi(t) \,dt
\Rightarrow
N\biggl(0,\int_{0}^1\varpi^2(t) \,dt\biggr),
\end{equation}
where $\varpi(t)=\pi(t)\sigma(t)/f(t,Q(t))$.
\end{lemma}
\begin{pf}
By Lemma \ref{lem::approximation}, we have
%
%
\begin{eqnarray}\label{eq::lem32}
&&{\sup_{t\in\mathcal{T}^*_n}}|f(t,Q(t))\mathfrak{Y}_n(t)-
\mathcal{S}_n(t)|\nonumber\\[-8pt]\\[-8pt]
&&\qquad=O_{\mathbb{P}} \biggl( {{\pi_n^{1/2}\log n+b_n^{3/4}} \over\sqrt{nb_n}}+
b_n\pi_n+ \pi_n^{2} \biggr).\nonumber
\end{eqnarray}

On the other hand, arguments similar to those in the proof of
(\ref{eq::thm113}) and (\ref{eq::thm114}) lead to
%
%
\begin{equation}\label{eq::lem35}
{\sup_{t\in\mathcal{T}^*_n}}|\mathcal{S}_n(t)-\sigma(t)\mathfrak{X}_n(t)|
=o_\mathbb{P}\biggl(\frac{n^{1/4}\log^2 n}{nb_n}+\frac{b_n\log n}{\sqrt
{nb_n}} \biggr),
\end{equation}
where
%
%
\begin{equation}\label{eq::lem36}
\mathfrak{X}_n(t)=\sum_{i=1}^nV_iK^*_{b_n}(t_i-t)/(nb_n).
\end{equation}
Recall $(V_i)_1^n$ are i.i.d. standard normal random variables defined
in (\ref{eq::thm110}). Hence, by (\ref{eq::lem32}), (\ref{eq::lem35})
and the bandwidth conditions of Theorem \ref{thm::L2test}, it follows that
%
%
\begin{equation}\label{eq::lem37}
\sup_{t\in\mathcal{T}^*_n} \biggl|\mathfrak{Y}_n(t)- \frac{\sigma
(t)}{f(t,Q(t))}\mathfrak{X}_n(t) \biggr|=o_\mathbb{P}(n^{-1/2}).
\end{equation}
Furthermore, it is easy to obtain that
%
%
\begin{equation}\label{eq::lem38}
\sqrt{n}\int_{\mathcal{T}_n^*}\mathfrak{X}_n(t)\varpi(t)
\,dt\Rightarrow
N\biggl(0,\int_{0}^1\varpi^2(t) \,dt\biggr).
\end{equation}
By (\ref{eq::lem37}) and (\ref{eq::lem38}), this lemma follows.
\end{pf}
\begin{lemma}\label{lem::quadratic}
Recall (\ref{L2normality}) for the definition of $\pi^*(t)$ and
(\ref{eq::lem36}) for the definition of $\mathfrak{X}_n(t)$. We have
\begin{eqnarray*}
&&
n\sqrt{b_n}\int_{\mathcal{T}_n^*}\mathfrak{X}_n(t)^2\pi^*(t)
\,dt-\frac{1}{\sqrt{b_n}}K^*\star K^*(0)\int_{0}^1\pi^*(t) \,dt\\
&&\qquad\Rightarrow N \biggl(0,2\int_\mathbb{R}[K^*\star K^*(t)]^2 \,dt\int
_{0}^1\pi
^*(t)^2 \,dt \biggr).
\end{eqnarray*}
\end{lemma}
\begin{pf}
Define
\begin{eqnarray*}
I_n&=&\sum_{i=1}^nR_i^2\int_{\mathcal{T}_n^*}[K^*_{b_n}(t_i-t)]^2\pi
^*(t) \,dt/(nb_n^{3/2}),\\
\mathit{I I}_n&=&\sum_{1\le i\neq j\le n}R_iR_j\int_{\mathcal
{T}_n^*}K^*_{b_n}(t_i-t)K^*_{b_n}(t_j-t)\pi^*(t) \,dt/(nb_n^{3/2}).
\end{eqnarray*}
Then by the central limit theorem for $I_n$, it follows that
%
%
\begin{eqnarray}\label{eq::lem41}
I_n&=&\frac{1}{\sqrt{b_n}}K^*\star K^*(0)\int_{0}^1\pi^*(t)
\,dt+O_\mathbb{P}
\biggl(\frac{1}{nb_n^{3/2}}+\frac{1}{\sqrt{nb_n}} \biggr)\nonumber\\[-8pt]\\[-8pt]
&=&\frac{1}{\sqrt{b_n}}K^*\star K^*(0)\int_{0}^1\pi^*(t)
\,dt+o_\mathbb{P}(1).\nonumber
\end{eqnarray}
On the other hand, by Theorem 2.1 of de Jone (\citeyear{deJone87}),
elementary but
tedious calculations show that
%
%
\begin{equation}\label{eq::lem42}
\mathit{I I}_n \Rightarrow N \biggl(0,2\int_\mathbb{R}[K^*\star
K^*(t)]^2 \,dt\int_{0}^1\pi^*(t)^2 \,dt \biggr).
\end{equation}
Since $n\sqrt{b_n}\int_{\mathcal{T}_n^*}\mathfrak{X}_n(t)^2\pi^*(t)
\,dt=I_n+\mathit{I I}_n$, by (\ref{eq::lem41}) and (\ref{eq::lem42}), the lemma follows.
\end{pf}
\begin{proposition}\label{prop::quad_clt}
Let $\tilde{T}_n=\int_{\mathcal{T}_n^*}\mathfrak{Y}^2_n(t)\pi(t) \,dt$.
Then under the conditions of Theorem \ref{thm::L2test}, we have
%
%
\begin{eqnarray}\label{eq::prop11}
&&n\sqrt{b_n}\tilde{T}_n-\frac{1}{\sqrt{b_n}}K^*\star K^*(0)\int
_{0}^1\pi^*(t) \,dt\nonumber\\[-8pt]\\[-8pt]
&&\qquad\Rightarrow N \biggl(0,2\int_\mathbb{R}[K^*\star K^*(t)]^2 \,dt\int
_{0}^1\pi
^*(t)^2 \,dt \biggr).\nonumber
\end{eqnarray}
\end{proposition}
\begin{pf}
By (\ref{eq::lem32}) and (\ref{eq::lem35}), we have
%
%
\begin{equation}\label{eq::prop12}
\sup_{t\in\mathcal{T}^*_n}\biggl|\mathfrak{Y}_n(t)- \frac{\sigma
(t)}{f(t,Q(t))}\mathfrak{X}_n(t)\biggr|=O_\mathbb{P}(\nu_n),
\end{equation}
where $\nu_n={{\pi_n^{1/2}\log n+b_n^{3/4}} \over\sqrt{nb_n}}+
b_n\pi_n+ \pi_n^{2}+\frac{n^{1/4}\log^2 n}{nb_n}$.

On the other hand, it is easy to show that
%
%
\begin{equation}\label{eq::prop14}
{\sup_{t\in\mathcal{T}^*_n}|\mathfrak{X}_n(t)}|=O_\mathbb{P}\biggl(\frac
{\log
^{1/2}n}{\sqrt{nb_n}}\biggr).
\end{equation}
Therefore, (\ref{eq::prop12}) and (\ref{eq::prop14}) imply that
%
%
\begin{equation}\label{eq::prop13}
{\sup_{t\in\mathcal{T}^*_n}}|\mathfrak{Y}_n(t)|=O_\mathbb{P}\biggl(\frac
{\log
^{1/2}n}{\sqrt{nb_n}}\biggr).
\end{equation}

Hence, by (\ref{eq::prop12}), (\ref{eq::prop14}), (\ref{eq::prop13})
and bandwidth conditions of Theorem \ref{thm::L2test},
%
%
\begin{eqnarray}\label{eq::prop15}
&&
\biggl|\tilde{T}_n-\int_{\mathcal{T}_n^*}[\mathfrak{X}_n(t)]^2\pi^*(t)\,dt\biggr|\nonumber\\
&&\qquad= \biggl|\int_{\mathcal{T}_n^*}\Delta_n(t) \biggl(\mathfrak{Y}_n(t)+
\frac{\sigma(t)}{f(t,Q(t))}\mathfrak{X}_n(t) \biggr)\pi(t) \,dt
\biggr|\nonumber\\[-8pt]\\[-8pt]
&&\qquad\le {C\sup_{t\in\mathcal{T}^*_n}}|\Delta_n(t)|\sup_{t\in\mathcal
{T}^*_n} \biggl[|\mathfrak{Y}_n(t)|+ \biggl|\frac{\sigma
(t)}{f(t,Q(t))}\mathfrak{X}_n(t) \biggr| \biggr]\nonumber\\
&&\qquad=O_\mathbb{P}\biggl(\nu_n \biggl(\frac{\log^{1/2}n}{\sqrt{nb_n}} \biggr)
\biggr)=o_\mathbb{P}
\biggl(\frac{1}{n\sqrt{b_n}} \biggr),\nonumber
\end{eqnarray}
where $\Delta_n(t)=\mathfrak{Y}_n(t)- \frac{\sigma
(t)}{f(t,Q(t))}\mathfrak{X}_n(t)$. Therefore, by Lemma \ref
{lem::quadratic}, we illustrate that Proposition \ref{prop::quad_clt} holds.
\end{pf}
\begin{pf*}{Proof of Theorem \ref{thm::L2test}}
First, note
\[
\tilde{Q}(t)-Q^o(t)=\mathfrak{Y}_n(t)+\mathbb{E}\mathcal
{S}^o_n(t)/f(t,Q(t))+\varrho_n\bigl(\eta(t)+o(1)\bigr).
\]
Therefore,
%
%
\begin{eqnarray}\label{eq::thm31}
T_n^*&=&\tilde{T}_n+\int_{\mathcal{T}_n^*}\bigl[\mathbb{E}\mathcal
{S}^o_n(t)/f(t,Q(t))+\varrho_n\bigl(\eta(t)+o(1)\bigr)\bigr]^2\pi(t) \,dt\nonumber\\
&&{} + 2\int_{\mathcal{T}_n^*}\mathfrak{Y}_n(t)\bigl[\mathbb{E}\mathcal
{S}^o_n(t)/f(t,Q(t))+\varrho_n\bigl(\eta(t)+o(1)\bigr)\bigr]\pi(t) \,dt\\
:\!&=&\tilde{T}_n+I^*_n+\mathit{I I}_n^*.\nonumber
\end{eqnarray}
Since $Q(t)\in\mathcal{C}^2[0,1]$, simple calculations show that
%
%
\begin{equation}\label{eq::thm33}
{\sup_{t\in\mathcal{T}_n^*}}|\mathbb{E}\mathcal
{S}^o_n(t)|=o(b_n^2)+O \biggl(\frac
{1}{nb_n} \biggr).
\end{equation}
By the bandwidth conditions of Theorem \ref{thm::L2test}, it follows that
%
%
\begin{equation}\label{eq::thm34}
I_n^*-\varrho_n^2\int_{\mathcal{T}_n^*}\eta^2(t)\pi(t)
\,dt=o(b^4_n+b_n^2\varrho_n)=o(\varrho_n^2).
\end{equation}
By Lemma \ref{lem::linearform},
%
%
\begin{eqnarray}\label{eq::thm32}
\varrho_n\int_{\mathcal{T}_n^*}\mathfrak{Y}_n(t)[\eta(t)+o(1)]\pi(t)
\,dt&=&O_\mathbb{P}(\varrho_nn^{-1/2})\nonumber\\[-8pt]\\[-8pt]
&=&o_\mathbb{P}(n^{-1}(b_n)^{-1/2}).\nonumber
\end{eqnarray}
Similarly,
%
%
\begin{equation}\label{eq::thm35}
\int_{\mathcal{T}_n^*}\mathfrak{Y}_n(t)[\mathbb{E}\mathcal
{S}^o_n(t)/f(t,Q(t))]\pi(t) \,dt=o_\mathbb{P}(n^{-1}(b_n)^{-1/2}).
\end{equation}
Hence, Proposition \ref{prop::quad_clt}, (\ref{eq::thm31}), (\ref
{eq::thm34}), (\ref{eq::thm32}) and (\ref{eq::thm35}) imply that Theorem
\ref{thm::L2test} holds.
\end{pf*}
\begin{remark}\label{rem::bias}
Under the null hypothesis $Q(t)=Q^o(t)$, we see from the above
proof that the bias $\mathcal{B}_n(t)$ influences the asymptotic
behavior of $T_n^*$ through two terms $\int_{\mathcal{T}_n^*}\mathcal
{B}^2_n(t)\pi(t) \,dt$ and $\int_{\mathcal{T}_n^*}\mathfrak
{Y}_n(t)\mathcal{B}_n(t)\pi(t) \,dt$, where $\mathcal{B}_n(t)=\mathbb{E}
\mathcal
{S}^o(t)/f(t,Q(t))$.
The jackknife bias reduction technique reduces those two effects to
second order. However, it can be shown that if the original estimate
$\hat{Q}(t)$ is used in the $\mathcal{L}^2$ test, then the first term
is not negligible under the optimal bandwidth $b_n=O(n^{-2/9})$,
which complicates the asymptotic analysis and reduces the precision
of the test.
\end{remark}
\begin{pf*}{Proof of Proposition \ref{prop::para_quantile}}
Proposition \ref{prop::para_quantile} follows from Lemmas \ref{lem::6}
and \ref{lem::5} below.
\end{pf*}
\begin{lemma}\label{lem::6}
Let $L_n=\sum_{i=1}^n\psi(e_i)\mathbf{g}(t_i)$. Then under the conditions
of Proposition \ref{prop::para_quantile}, we have $L_n=O_\mathbb
{P}(\sqrt{n})$.
\end{lemma}
\begin{pf}
Let $L_{n,k}=\sum_{i=1}^n\mathcal{P}_{i-k}\psi(e_i)\mathbf
{g}(t_i)$. Recall the
operator $\mathcal{P}_{k}$ is defined in Theorem \ref{thm::sip}. Then the
summands of $L_{n,k}$ are martingale differences. By the orthogonality,
we have for $k\ge1$
\begin{eqnarray*}\label{eq::lem61}
\|L_{n,k}\|^2=\sum_{i=1}^n\|\mathcal{P}_{i-k}\psi(e_i)\mathbf
{g}(t_i)\|
^2&=&\sum
_{i=1}^n|\mathbf{g}(t_i)|^2\|\mathcal{P}_{i-k}\psi(e_i)\|^2\nonumber
\\
&\le& Cn\delta^2_F(k-1,2).
\end{eqnarray*}
Similar arguments also imply $\|L_{n,k}\|^2\le Cn$. Therefore,
\[
\|L_n\|=\sum_{k=0}^{\infty}\|L_{n,k}\|\le C\sqrt{n}\Biggl(1+\sum
_{k=1}^{\infty
}\delta_F(k-1,2)\Biggr)\le C\sqrt{n}.
\]
Hence, Lemma \ref{lem::6} follows.
\end{pf}
\begin{lemma}\label{lem::5}
Under the conditions of Proposition \ref{prop::para_quantile}, we
have
%
%
\begin{equation}\label{eq::lem51}
\mathcal{H}_{\alpha}\tilde{\theta}_{\alpha}-\sum_{i=1}^n\psi
_{\alpha
}(e_i)\mathbf{g}(t_i)/\sqrt{n}=o_p(1),
\end{equation}
where $\tilde{\theta}_{\alpha}=\sqrt{n}(\hat{\theta}_{\alpha
}-\theta
_0)$, $e_i=X_i-Q_{\alpha}(t_i)$ and
\[
\mathcal{H}_{\alpha}=\int_{0}^1\mathbf{g}(t)\mathbf{g}^{\top
}(t)f(t,Q_{\alpha
}(t)) \,dt.
\]
\end{lemma}
\begin{pf}
We shall omit the subscript $\alpha$ in the proof. Let $\mathbf
{g}_n(t)=\mathbf{g}(t)/\sqrt{n}$. By arguments similar to those of Lemma
3 of Wu (\citeyear{W07a}), we have for any fixed constant $c$ and fixed
$\theta\le c$,
%
%
\begin{equation}\label{eq::lem52}
\operatorname{var}\Biggl(\sum_{i=1}^n\eta_i(\theta)\Biggr)=o(1),
\end{equation}
where $\eta_i(\theta)=\rho(e_i-\mathbf{g}^{\top}_n(t_i)\theta
)-\rho
(e_i)+\mathbf{g}^{\top}_n(t_i)\theta\psi(e_i)$. On the other hand,
elementary calculations based on the Taylor expansion show that
%
%
\begin{eqnarray}\label{eq::lem53}
\sum_{i=1}^n\mathbb{E}[\eta_i(\theta)]&=&\sum_{i=1}^n \biggl[\frac
{f(t_i,Q(t_i))}{2}|\mathbf{g}^{\top}_n(t_i)\theta|^2+o(|\mathbf
{g}^{\top
}_n(t_i)\theta|^2) \biggr]\nonumber\\[-8pt]\\[-8pt]
&=&\frac{\theta^{\top}\mathcal{H}\theta}{2}+o(1).\nonumber
\end{eqnarray}
From equations (\ref{eq::lem52}), (\ref{eq::lem53}) and the
convexity lemma [Pollard (\citeyear{P91}), page 187], we obtain
%
%
\begin{equation}\label{eq::lem54}\qquad
\sup_{\theta\le c} \Biggl|\sum_{i=1}^n \bigl[\rho\bigl(e_i-\mathbf{g}^{\top
}_n(t_i)\theta\bigr)-\rho(e_i)+\mathbf{g}^{\top}_n(t_i)\theta\psi(e_i)
\bigr]-\frac{\theta^{\top}\mathcal{H}\theta}{2} \Biggr|=o_\mathbb{P}(1).
\end{equation}
Now a standard argument using properties of convex functions will
lead to (\ref{eq::lem51}). See, for example, the proofs of Theorems
2.2 and 2.4 in Bai, Rao and Wu (\citeyear{BRW92}). Details are omitted.
\end{pf}

\section*{Acknowledgments}
We are grateful to two referees, the
Associate Editor and the editor for their many helpful comments.

\printaddresses

\end{document}